\documentclass[runningheads]{llncs}
\usepackage{authblk}
\usepackage[utf8]{inputenc}
\usepackage{graphicx}
\usepackage{subcaption}
\usepackage{amsmath,amsfonts,mathtools}
\usepackage{derivative}
\usepackage{svg}
\usepackage{algorithm}
\usepackage{algcompatible}
\usepackage{algpseudocode}
\algblockdefx{ParFor}{EndParFor}[1]%
{\textbf{parfor }#1 \textbf{ do}}%
{\textbf{end parfor}}
\DeclareMathOperator{\diag}{diag}

\usepackage{subdepth} 
\usepackage{mathtools}
\usepackage{multirow}
\usepackage{hhline}
\usepackage{colortbl}
\usepackage{diagbox} 
\usepackage{makecell} 
\usepackage[misc,geometry]{ifsym}
\newcommand{\B}[1]{\mbox{\boldmath $#1$}}

\begin{document}

\title{Experimental Study of a Parallel Iterative Solver for Markov Chain Modeling}

\author{V. Besozzi \inst{1}\and
M. Della Bartola \inst{1}\and 
L. Gemignani \inst{1}\orcidID{0000-0001-8000-4906}$\spadesuit$ }

\institute{Dipartimento di Informatica, Universit\`a di Pisa, 56127, Pisa, Italy \\
$\spadesuit$ Corresponding author: \email{luca.gemignani@unipi.it}}


\maketitle
\begin{abstract}
This paper presents the results of a preliminary experimental investigation of the  performance of a stationary iterative method based on a block staircase splitting  for solving singular systems of linear equations arising in Markov chain modelling.  From the experiments presented, we can deduce  that the 
method is well suited for solving block banded or more generally localized systems  in a parallel computing environment.  The parallel implementation has been benchmarked using several Markovian models.

\end{abstract}

\section{Introduction}

The solving of linear algebraic systems lies at the core of many scientific and engineering simulations. Discrete-state  models are  widely employed for modeling and analysis of large networks and systems such as communication networks, allocation schemes, computer
systems and population processes. If the future evolution of the system depends
 only on the current state of the system and not on the
 past history, the system may be represented by a
 Markov  chain.   For a homogeneous, irreducible, continuous time
 Markov chain  with $N$ states, the long-term behaviour of the system is determined by the  stationary probability vector $\B \pi \in \mathbb R^N$ such that 
 \begin{equation}\label{maineq}
 Q^T\B \pi=\B 0, \quad  \B \pi\geq \B 0, \quad \B e^T\B \pi=1, 
 \end{equation}
 where $Q\in \mathbb R^{N\times N}$  is the transition rate matrix, or the infinitesimal generator of the Markov chain,  and $\B e=\left[1, \ldots, 1\right]^T$.   Since $Q$ is irreducible, by the Perron-Frobenius Theorem \cite{Mey} we  find that $Q$ has rank $N-1$ and, therefore, $\B \pi$ spans the kernel of $Q^T$.  The computation of $\B \pi$ amounts to solve the homogeneous linear system  \eqref{maineq}. A review of numerical methods for solving \eqref{maineq} can be found in \cite{SWJ,PSS}. Active research  in this area is focused on the development of techniques, methods and data structures, which minimize  the computational (space and
time) requirements for solving the linear system \eqref{maineq} when $Q$ is large and sparse. One of such techniques is parallelization.

 Iterative methods  are generally preferred for  solving  large linear systems of equations   because they are insensitive to fill-in and  accuracy issues \cite{PSS}. Stationary iterative methods like Gauss-Seidel (GS), Jacobi,
 and Successive Over-Relaxation (SOR)  are  interesting on their own and have further applications as preconditioners for   projection methods  like CG and GMRES.  Experimental studies  
 demonstrated that for Markov chain problems  \eqref{maineq} block methods based on matrix splittings such as block Jacobi and block Gauss–Seidel give better convergence than  other projection methods  (see \cite{Tou} and the references given therein).

 Among classical iterative methods, the Gauss-Seidel method has several interesting  features.
It is a classical result that on a nonsingular M-matrix the Gauss-Seidel method converges faster than the Jacobi method~\cite[Corollary~5.22]{Plemmons}. Moreover it can be implemented just using one iteration vector which is an important feature for  huge systems. 
The SOR method with the optimal relaxation parameter can be better yet, but, however, choosing an optimal SOR relaxation parameter is  difficult for many problems.
Therefore, the Gauss--Seidel method is very attractive in practice and it is  also used  as preconditioner in combination with other iterative schemes.  A classical example is the
multigrid method for partial differential equations, where using Gauss--Seidel or SOR as a smoother typically yields good convergence properties \cite{DBLP}.
 
Parallel implementations of Gauss–Seidel method have been designed for  certain regular problems, for example, the solution of Laplace’s equations by finite differences,
by relying upon  red-black coloring or more generally multi-coloring  schemes to provide some parallelism \cite{OV}.
In most cases, constructing efficient parallel true  Gauss–Seidel algorithms is challenging  and ProcessorBlock (or localized) Gauss–Seidel is often used \cite{SY}. Recent examples with applications to Markov chain modeling are the methods proposed in \cite{BB}  and \cite{AMKS}.
Here, each processor performs Gauss–Seidel as a subdomain solver for a block Jacobi method. While Processor Block Gauss–Seidel  methods are  easy to parallelize, the overall convergence  of the resulting iterative scheme  can suffer.

In order to cope with the parallelization of Gauss-Seidel type methods while retaining the same convergence rate in  \cite{S1}  staircase  splittings were introduced by  proving  that
for consistently ordered matrices \cite{Saad}  the iterative scheme based on such partitionings splits  into  independent computations and  at the same time exhibits the same convergence rate as the
classical Gauss-Seidel iteration.  A specialization  of this result  for  block tridiagonal matrices  had  already appeared in \cite{AM}.   More recently, in \cite{GP}   the computational interest of  staircase splittings has been broadened  by showing  that  for a nonsingular  M-matrix $A$ in block  lower Hessenberg  form   the asymptotic rate of convergence of the block staircase method is better than the asymptotic rate of convergence of the block Gauss-Seidel  method  applied to $A$.  A further extension with applications  to accelerating certain fixed point iterations for Markov chain modeling is 
given in \cite{GM_NUMA}. These results are  quite surprising since the matrix $M$ in the block staircase partitioning of $A=M-N$ is much more sparse than the corresponding block lower triangular matrix  $M$ of the Gauss-Seidel splitting. Moreover,  our  experimental evidence indicates that the block staircase partitioning generally works quite well   when compared to block Gauss-Seidel for   block  banded or more generally localized matrices \cite{Benzi} with  entries decaying away from the  main diagonals. 

The contribution of this paper  is twofold.    The matrix $A=Q$  in \eqref{maineq} is singular and   the  comparison theorems proved in \cite{GP} do not extend to the singular case   while classical results for singular systems \cite{MS0,MS1} do not apply to our methods.  The first aim is to  gain an understanding  of how (block) staircase  and (block) Gauss-Seidel type methods  compare  when applied for solving  large and sparse Markov Chain problems.   In particular, we are interested in the case where $A$ is banded or localized around the main diagonals.  The second goal is to perform this comparison in a parallel computing environment.  To do this we have implemented  a block staircase iterative solver  for the parallel  computation of the vector $\B \pi$.  The properties of this  method are examined experimentally.  In particular,    numerical experiments are performed to compare our  method  with an implementation of the composite solver proposed in \cite{BB}  in terms of traditional efficiency measures for parallel algorithms.   A discussion of the results is presented  together with some conclusions and  insights for future work. 

\section{Mathematical Background}

Let $P=(p_{i,j})\in \mathbb R^{N\times N}$  be a transition probability matrix of a  homogeneous  ergodic Markov Chain with $N$ states.   Then $P$ is  irreducible and row-stochastic, that is, $p_{i,j}\geq 0$, $1\leq i,j\leq N$, and $P\B e=\B e$ with $\B e^T=\left[1, \ldots, 1\right]$.  The matrix $Q^T=I_N-P^T$ is a singular M-matrix.  Observe that $ \B e^TQ^T=\B 0^T$.  Since $Q^T$ is also irreducible, by the Perron-Frobenius Theorem  it follows that  the  kernel of $Q^T$ is spanned by a  vector $\B \pi$  such that $\B \pi>\B 0$ and $\B e^T \B \pi =1$.  This vector is called the stationary probability distribution vector of the Markov Chain. 

The computation of $\B \pi$ amounts to solve the  homogeneous linear system $Q^T \B \pi=\B 0$ under the normalization $\B e^T \B \pi=1$.  Iterative methods  based on the power iteration can  be used \cite{PSS}. The computational efficiency and the convergence properties of these algorithms can benefit of a  block partitioning of the matrix $Q^T$.  Let us assume that 
\[
Q^T=\left[\begin{array}{ccc}Q_{1,1} & \ldots & Q_{1,n}\\\vdots   & \vdots & \vdots \\Q_{n,1} & \ldots & Q_{n,n}\end{array}\right], 
\]
where $Q_{i,j}\in \mathbb R^{n_i\times n_j}$, $1\leq i,j\leq n$, $\sum_{i=1}^n n_i=N$. 

A regular splitting of the matrix $Q^T$ is a partitioning $Q^T=M-N$ with 
$M^{-1}\geq 0$ and $N\geq 0$. 
Since $Q^T\pi=\B 0$ we find that $M\B \pi=N \B\pi$
which gives  $M^{-1} N \B \pi=\B \pi$. It is well known that the spectral radius  of  $M^{-1} N$ is equal to 1 and $\lambda=1$ is a simple eigenvalue of 
$M^{-1}N$ \cite{Sch}.  This not immediately implies that $\lambda=1$ is the dominant eigenvalue of $M^{-1}N$, that is,  that for the remaining eigenvalues $\lambda$ of $M^{-1}N$  it holds $|\lambda|<1$. 

\begin{example}\label{ex1} 
Let $Q^T=\left[\begin{array}{cccc} 1 & -1& 0& 0\\-1/2 & 1 & -1/2 & 0 \\ 0 & -1/2 & 1 & -1/2 \\ 0 & 0 & -1 & 1\end{array}\right]$.  The Jacobi splitting with $M=I_4$ is a regular splitting  but the iteration matrix $M^{-1} N$  has eigenvalues $\{-1, 1,-1/2, 1/2\}$.  The Gauss-Seidel splitting  
is a regular splitting and the   corresponding iteration matrix 
has eigenvalues $\{1, 1/4,0, 0\}$. 
\end{example}

By graph-theoretic arguments \cite{Sch} it follows that for a regular splitting the matrix $M^{-1}N$
is permutationally similar to a  block matrix $T=\left[\begin{array}{cc}0 & T_{1,2} \\ 0 & T_{2,2}
\end{array}\right]$ where $T_{2,2}$ is square, irreducible  and non-negative and  every row of the possibly nonempty matrix $T_{1,2}$  is nonzero.  A non-negative  square  matrix $A$  is primitive 
if there is $m\geq 1$ such that $A^m>0$.  By the Perron-Frobenius Theorem  we obtain that $\lambda=1$
is the dominant eigenvalue of $ M^{-1} N$ if $T_{2,2}$  is primitive.  Hereafter, this condition is 
always assumed.  Under this assumption  the classical power iteration is eligible for determining a numerical approximation of the vector $\B \pi$.

The method based on the (block) Jacobi splitting  is very  convenient  to vectorize and  to parallelize.  As shown in the simple example above it can suffer  from convergence  problems.  In this respect, the  (block) Gauss-Seidel iteration generally outperforms the Jacobi algorithm. 
ProcessorBlock (or localized) Gauss–Seidel schemes  provide a reliable compromise between parallelization and convergence issues.  One such   hybrid adaptation is described in \cite{BB}.  Suppose that  the matrix $Q^T$ is partitioned as 
\[
Q^T=\left[\begin{array}{c} {Q^{(1)}}^T\\\hline \vdots \\\hline {Q^{(p)}}^T\end{array}\right], \quad 
{Q^{(j)}}^T=\left[\begin{array}{ccc}Q_{m_j,1} & \ldots & Q_{m_j,n}\\\vdots   & \vdots & \vdots \\Q_{m_j+r_j-1,1} & \ldots & Q_{m_j+r_j-1,n}\end{array}\right], 
\]
with $m_1=1$, $m_j=\sum_{i=1}^{j-1}r_i+1$, $2\leq j\leq p$, $\sum_i^p r_i=n$. The  iterative 
scheme  in \cite{BB} exploits the regular splitting where 
\[
M=\diag\left[M_1, \ldots, M_p\right],\  M_j=\left[\begin{array}{ccc}Q_{m_j,m_j} \\
Q_{m_j+1,m_j}\!\! &\!\!Q_{m_j+1,m_j+1} \\\vdots \!\!& \!\! \ddots \\
Q_{m_{j+1}-1,m_j}\!\! & \!\!\ldots & \!\! Q_{m_{j+1}-1,m_{j+1}-1} \end{array}\right]
\]
and, hence, $M$ is a block diagonal matrix with  block lower triangular blocks.  The resulting scheme
proceeds as follows: 
\begin{equation}\label{mainsch}
\left\{ \begin{array}{ll} M\B x^{(k+1)}=N\B x^{(k)} \\
\B x^{(k+1)}=\frac{\B x^{(k+1)}}{\B e^T\B x^{(k+1)}} \end{array}\right., \quad k\geq 1. 
\end{equation}
If $p$ is the number of processors, then  {\bf {Algorithm 1}} is a possible implementation  of this scheme starting from the skeleton proposed in \cite{BB}. 

\begin{algorithm}
    \caption{This algorithm approximates the vector $\B \pi$ by means of the method in \cite{BB}}
    \label{algorithm1}
    \begin{algorithmic}[1] 
        \State Initialization
            \While{$tol\geq  err \ \& \ it\leq maxit $}
            \ParFor{$j=1, \ldots, p$} 
            \For{$k=m_j, \ldots, m_j+r_j-1$}
            \State $\mathcal J_k=\{s\in \mathbb N \colon m_j\leq s<k\}$
            \State $\B z_k\gets -Q_{k,k}^{-1}\left( \sum_{s\in \mathcal J_k} Q_{k,s}\B z_s +\sum_{s\in \{1,\ldots,n\}\setminus \mathcal J_k} Q_{k,s}\B x_s\right)$
            \EndFor
            \EndParFor
            \State $\B z \gets \frac{\B z}{\B e^T\B z}$
            \State $err\gets \parallel \B z-\B x\parallel_1$; $\B x\gets \B z$; $it \gets it+1$
            \EndWhile
            \State \textbf{return} $\B x$
    \end{algorithmic}
\end{algorithm}
A  different approach to parallelizing  stationary iterative solvers  was taken in \cite{S1}.  The approach is based on the  exploitation of  a suitable partitioning of the matrix $Q^T$ where  $M$ has a "zig-zag" pattern around the main diagonal referred  to  as a staircase profile.  More specifically,  we can associate  with $Q^T$ the stair  matrices $M_1\in \mathbb R^{N\times N}$ and $M_2 \in \mathbb R^{N\times N}$ defined by 
\[
M_1=\left[\begin{array}{cccccccccc}
    Q_{1,1} \\
    Q_{2,1} & Q_{2,2}& Q_{2,3}\\
    &   & Q_{3,3} \\
    & & Q_{4,3}& Q_{4,4}& Q_{4,5}\\
    &&&&  Q_{5,5}\\
    &&&&\times & \times  & \times  \\
    &&&&&&\times &\phantom{a}
    \\
    \\
  \end{array}\right]
\]
and 
\[
M_2=\left[\begin{array}{cccccccccc}
    Q_{1,1} & Q_{1,2}\\
     & Q_{2,2}\\
    &  Q_{3,2} & Q_{3,3} & Q_{3,4} \\
    & & & Q_{4,4}\\
    &&& Q_{5,4}&  Q_{5,5} & Q_{5,6}\\
    &&&&&\times   \\
    &&&&&\times & \times & \times\phantom{a}
    \\
    \\
  \end{array}\right]. 
\]
 These matrices  are called  stair matrices of type 1 and 2, respectively.
 Staircase splittings of the form $Q^T=M-N$ where $M$ is a stair matrix  have  two remarkable  features:
 \begin{enumerate}
     \item  The solution of a linear system $M\B x=\B b$   can be carried out in two parallel steps  since  all even and  all odd components of $\B x$ can be computed concurrently. 
     \item  In terms of convergence these splittings inherit some  advantages of the block  Gauss-Seidel method.  If $Q^T$ is block tridiagonal, then  it can be easily proved that the iteration matrices associated with   block  Gauss-Seidel and block staircase splittings  have the same eigenvalues and therefore the same convergence rate.  In \cite{GP} it is shown that the spectral radius  of the iteration matrix  generated  by the staircase splitting of an invertible M-matrix in  block lower Hessenberg form  is not greater than  the spectral radius of the corresponding  iteration matrix  in the block Gauss-Seidel method.  For singular matrices the convergence of the scheme \eqref{mainsch}  depends on the spectral gap between the dominant eigenvalue equal to  $1$ and the second eigenvalue 
     $\gamma=\max\{|\lambda| \colon \lambda \ {\rm {eigenvalue}} \ {\rm {of }} \  M^{-1}N \ {\rm {and }} \  \lambda \neq 1\}$. Experimentally   (see Example \ref{ex3} below) 
     the block staircase iteration  still performs similarly with the block Gauss-Seidel method when applied for solving linear systems with  singular M-matrices in block Hessenberg form.   The same behaviour is observed   for matrices that are localized around the main diagonals. Examples are the covariance matrices  with application to the spatial kriging problem  (compare with \cite{HSG}).  
     \end{enumerate}
\begin{example}\label{ex2}
For the matrix of Example \ref{ex1} the staircase splitting of the first type gives  an iteration matrix having the same eigenvalues  $\{1, 1/4,0, 0\}$ of the Gauss-Seidel  scheme. 
\end{example}
\begin{example}\label{ex3} 
We have performed several numerical experiments with  randomly generated singular M-matrices in banded block lower Hessenberg  form  having the profile depicted in Figure \ref{f1}. 

\medskip

  \begin{minipage}{\textwidth}
  \begin{minipage}[b]{0.49\textwidth}
    \centering
    \includegraphics[width=0.8\textwidth]{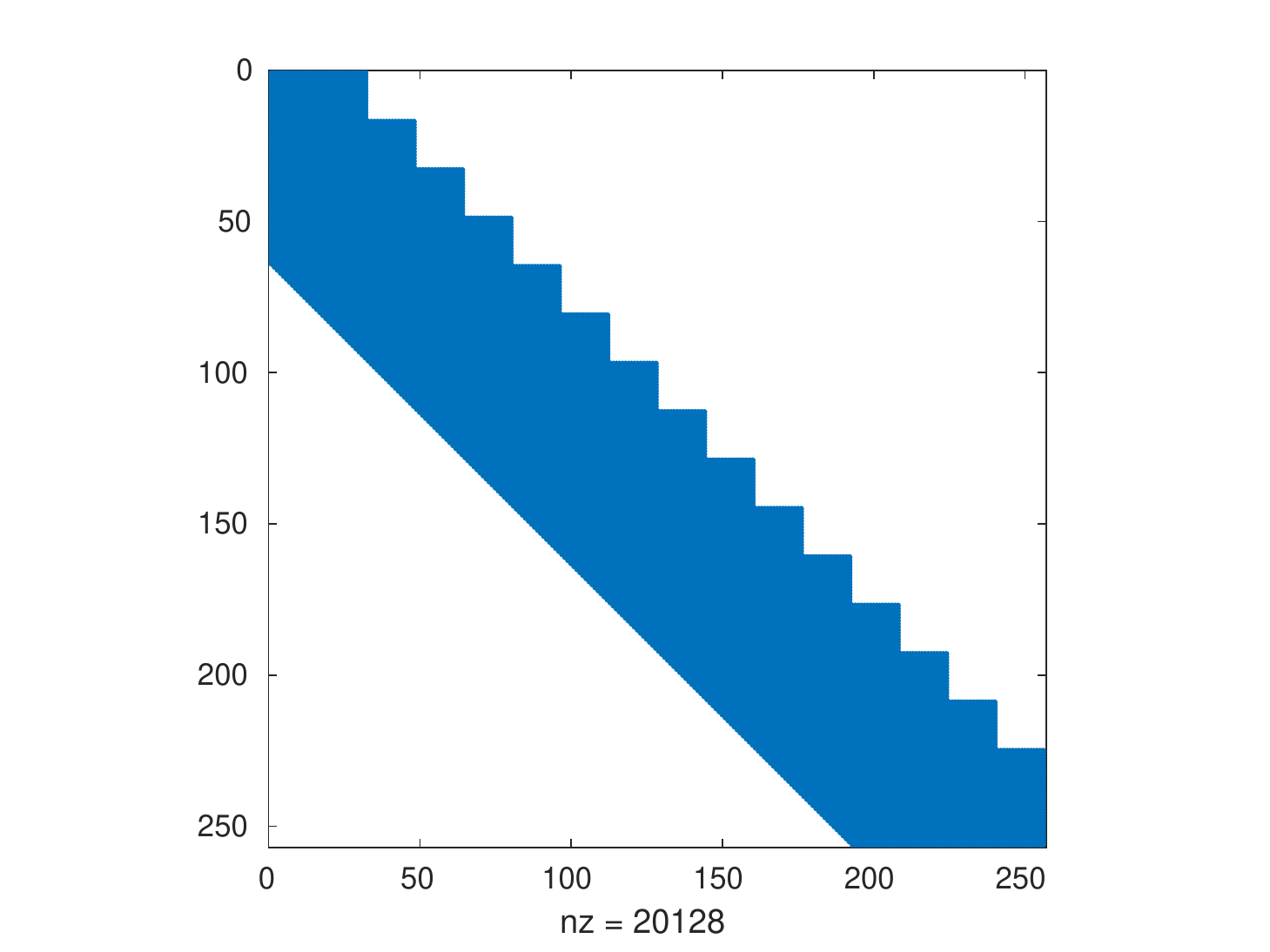}
    \captionof{figure}{Spy plot}
     \label{f1} 
  \end{minipage}
  \hfill
  \begin{minipage}[b]{0.49\textwidth}
    \centering
\begin{tabular}{|c|c||c|c|c|}
\hhline{--||---}
{\backslashbox{$kp$}{$k$}}
&  &  16 & 32& 64 \\
\hhline{==::===}
\multirow{2}{*}{16}
&$\rho_1$ &  &  & 1.6\\
\hhline{~-||---}
&$\rho_2$ &  &  & 1.5\\
\hhline{--||---}
\multirow{2}{*}{32}
&$\rho_1$ &  &  & 1.6\\
\hhline{~-||---}
&$\rho_2$ &  &  & 1.4\\
\hhline{--||---}
\multirow{2}{*}{64}
&$\rho_1$ & $\infty$ & 33.6 & 2.0\\
\hhline{~-||---}
&$\rho_2$ & 1.0 &1.0  & 1.4\\
\hhline{--||---}
\end{tabular}
\captionof{table}{Convergence comparisons}
\label{t1}
    \end{minipage}
  \end{minipage}
\medskip

The size of the blocks is $k$, the block size of the matrix is $n$ and  the lower bandwidth is $4k$.  We compare the performance of the block Jacobi ($BJ$), block Gauss-Seidel ($BGS$) and block staircase ($BS$) methods for different sizes of the block partitioning denoted as $kp$.  In  the following Table \ref{t1} we show the maximum value  of $\rho_1=\log(1/\gamma_{BGS})/\log(1/\gamma_{BJ})$ and of $\rho_2=\log(1/\gamma_{BGS})/\log(1/\gamma_{BS})$ --where $\gamma_M$ is the second eigenvalue of the iteration matrix generated by the method $M$-- over 1000  experiments with $k=64$ and $n=16$. The value $\infty$ indicates that in some  trials  the block Jacobi method does not converge due to the occurrence of two or more eigenvalues equal to 1 in magnitude. The results  demonstrate that the spectral gap of the iteration matrix in  the block staircase method remains close to  that one of BGS. 

\end{example}

The following {\bf Algorithm 2} provides an implementation of \eqref{mainsch} using  $M=M_1$.
 For $i=1, \ldots, n$ let 
 \[
 \mathcal J_i=\left\{\begin{array}{ll} \{i-1, i+1\}\cap \{1, \ldots, n\},   \ i  \ {\rm even}\\
 i,    \ i  \ {\rm odd.} \end{array}\right.
 \]

\begin{algorithm}
    \caption{This algorithm approximates the vector $\B \pi$ by means of the method \eqref{mainsch} with $M=M_1$}
    \label{algorithm2}
    \begin{algorithmic}[1] 
        \State Initialization
            \While{$tol\geq  err \ \& \ it\leq maxit $}
            \ParFor{$i=1, \ldots, n$} 
            \State $\B z_i\gets -\sum_{k\notin \mathcal J_i} Q_{i,k} \B x_k$
            \EndParFor
            \ParFor {$i=1,3, \ldots, n$} 
            \State $\B z_i\gets Q_{i,i}^{-1} \B z_i$
            \EndParFor
            \ParFor {$i=2,4, \ldots, n$} 
            \State $\B z_i\gets Q_{i,i}^{-1}\left(\B z_i-\sum_{k\in \mathcal J_i} Q_{i,k} \B z_k\right)$
            \EndParFor
            \State $\B z \gets \frac{\B z}{\B e^T\B z}$
            \State $err\gets \parallel \B z-\B x\parallel_1$; $\B x\gets \B z$; $it \gets it+1$
            \EndWhile
            \State \textbf{return} $\B x$
    \end{algorithmic}
\end{algorithm}
In the next section numerical experiments are  performed to compare  the performance of  {\bf Algorithm 1} and {\bf Algorithm 2} for solving singular systems  of linear equations  arising in Markov  chain modeling. 

\section{Numerical  Experiments}
The experiments have been run on a server with two Intel Xeon E5-2650v4 CPUs with 12 cores and 24 threads each, running at 2.20GHz. The parallel implementations of {\bf Algorithm 1} --referred to as {\tt JGS}  algorithm-- and {\bf Algorithm 2}  --referred to as {\tt STAIR1} or {\tt STAIR2} algorithm depending on the staircase splitting-- are based on OpenMP. Specifically, we have  used C++20 with the help of  Armadillo \cite{CC1,CC2} which also provides integration with LAPACK \cite{AndeBaiBiscBlacDemmDongDuCrGreeHammMcKeSore99}and OpenBLAS \cite{10.1145/355841.355847}.

\begin{figure}[ht] 
  \begin{subfigure}[b]{0.5\linewidth}
    \centering
    \includegraphics[width=0.75\linewidth]{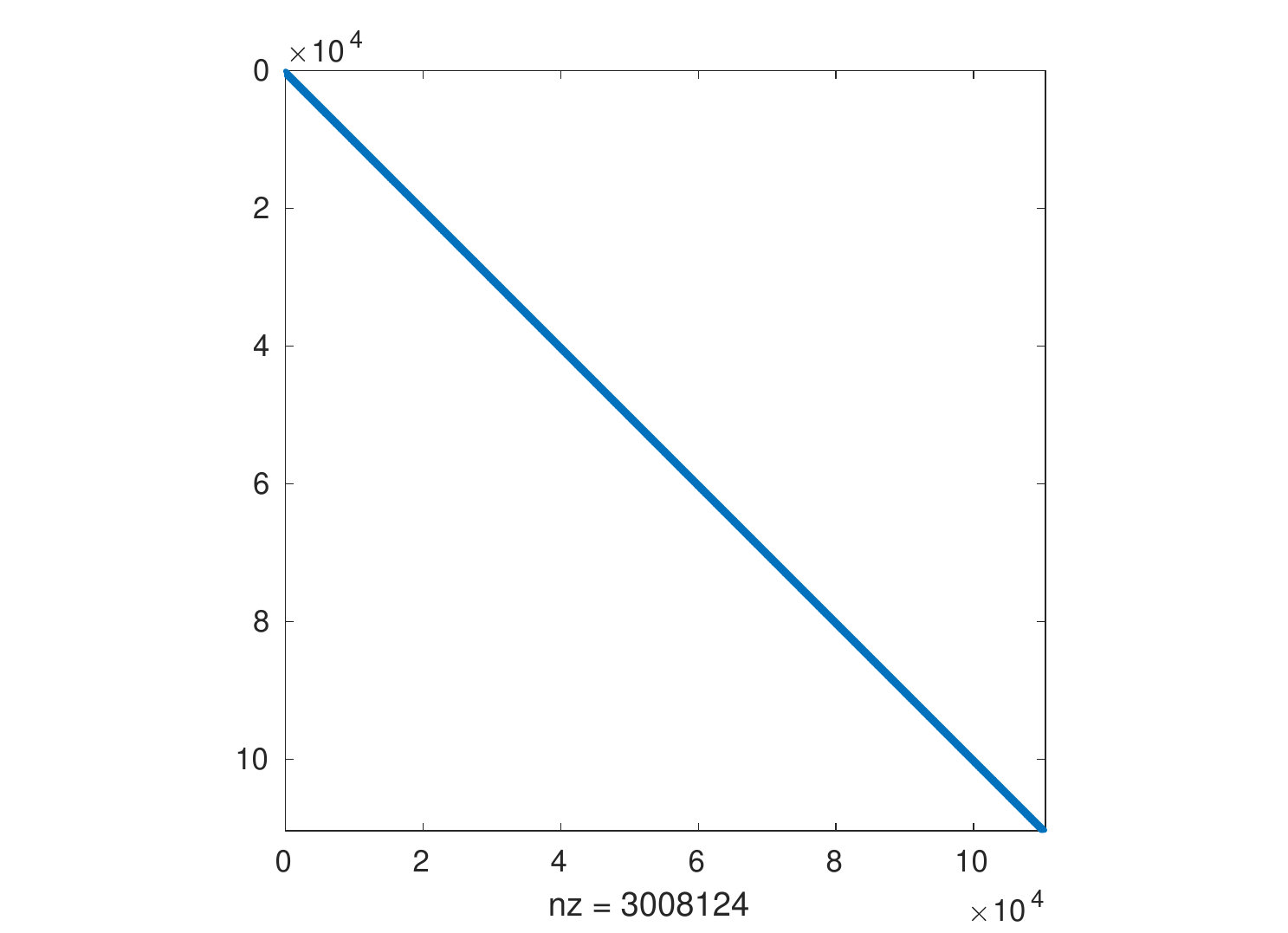} 
    \caption{Spy plot of the matrix in  \ref{it1}} 
    \label{fig7:a} 
    \vspace{4ex}
  \end{subfigure}
  \begin{subfigure}[b]{0.5\linewidth}
    \centering
    \includegraphics[width=0.75\linewidth]{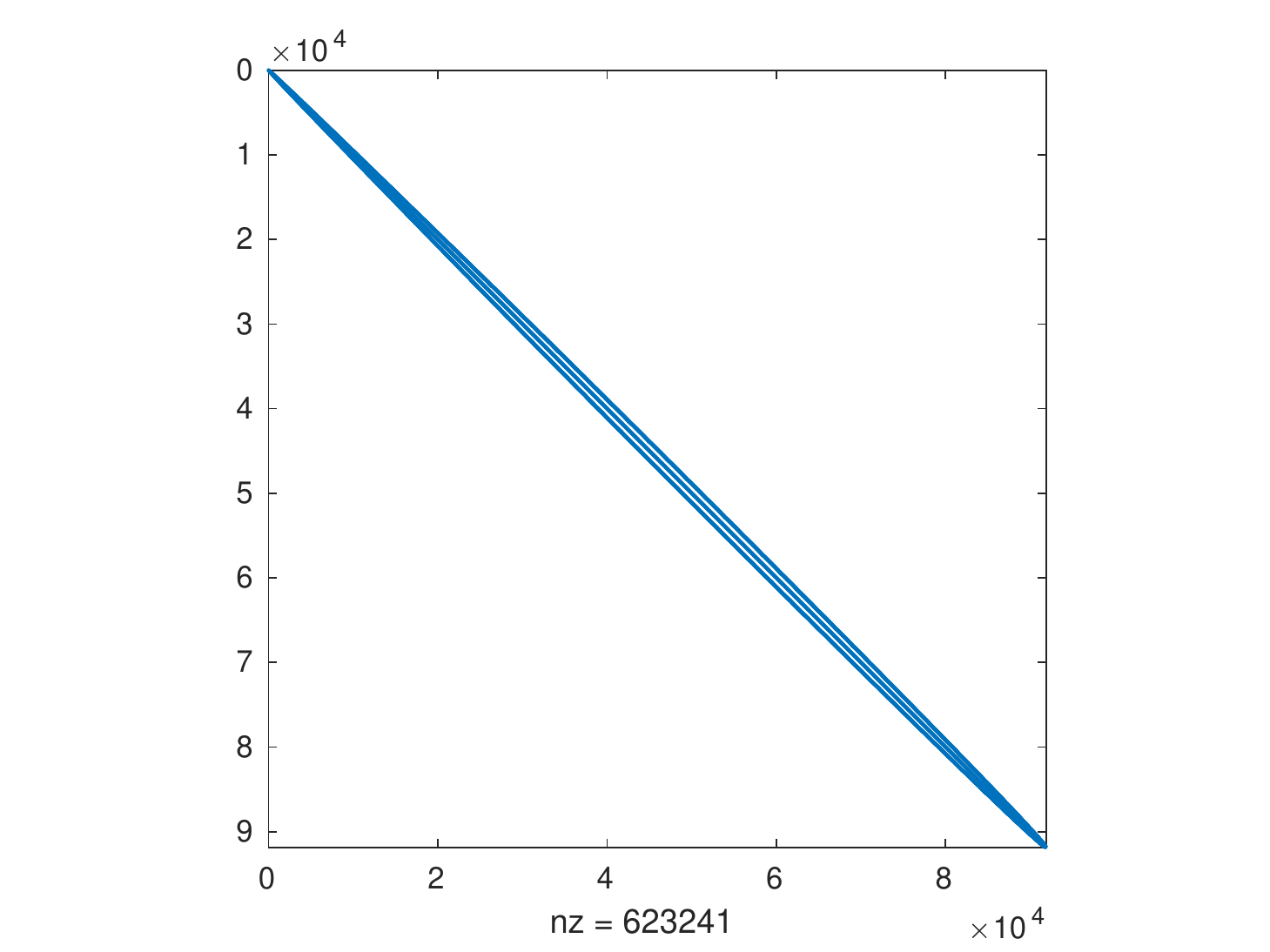} 
    \caption{ Spy plot of the matrix in  \ref{it2}}
    \label{fig7:b} 
    \vspace{4ex}
  \end{subfigure} 
  \begin{subfigure}[b]{0.5\linewidth}
    \centering
    \includegraphics[width=0.75\linewidth]{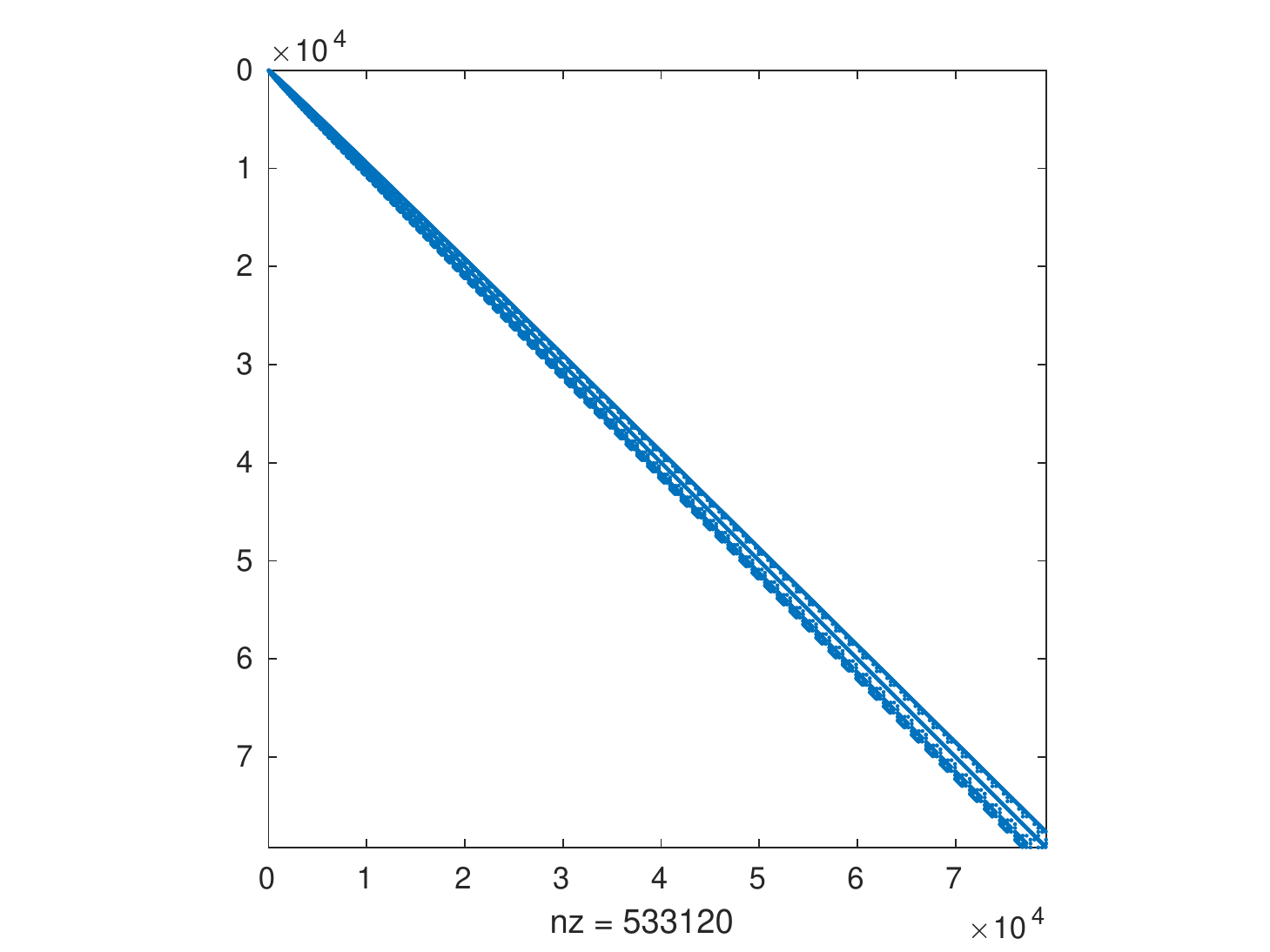} 
      \caption{Spy plot of the matrix in  \ref{it3}} 
    \label{fig7:c} 
  \end{subfigure}
  \begin{subfigure}[b]{0.5\linewidth}
    \centering
    \includegraphics[width=0.75\linewidth]{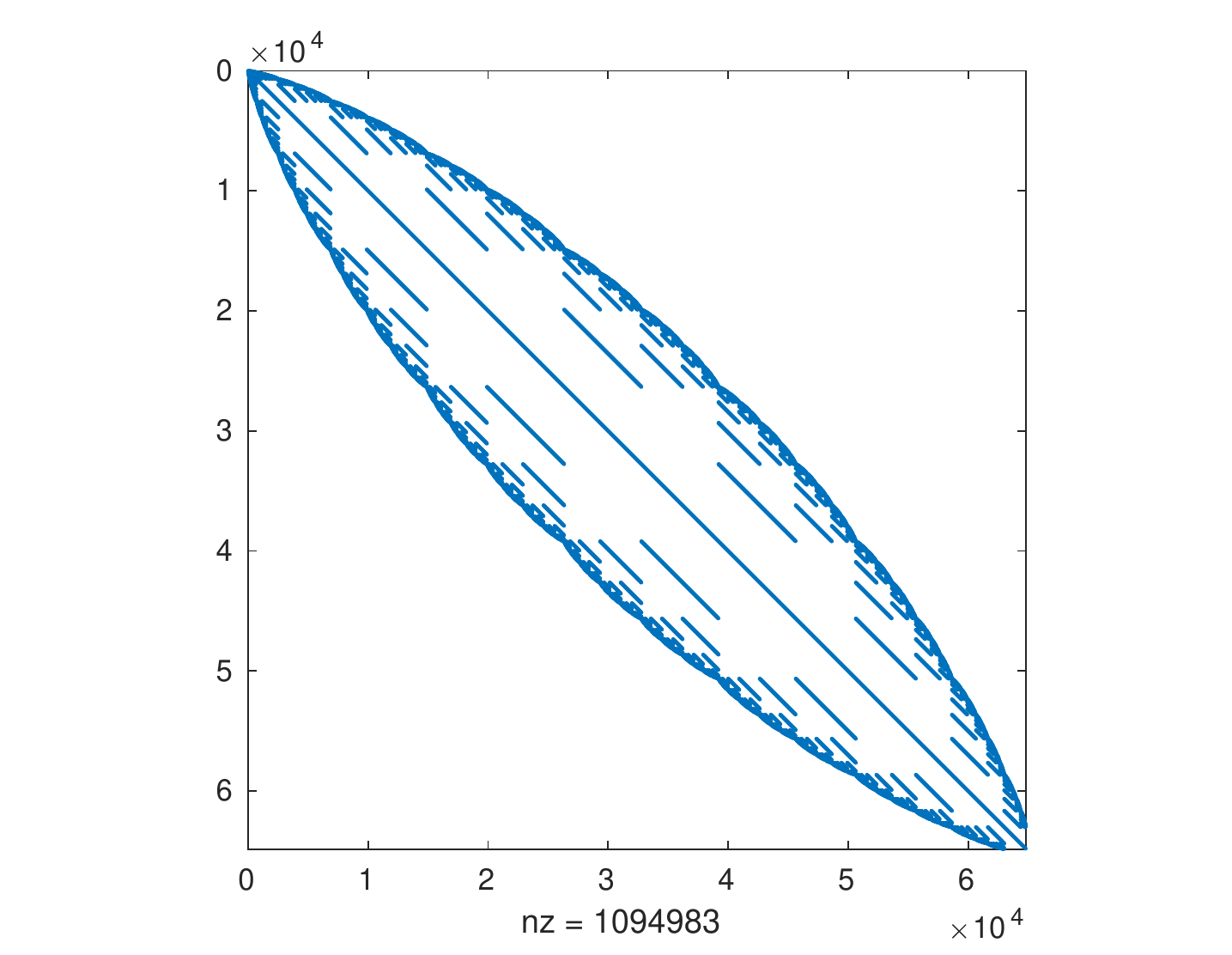} 
    \caption{ Spy plot of the matrix in  \ref{it4}}
    \label{fig7:d} 
  \end{subfigure} 
  \caption{Illustration of   the  sparsity pattern of the matrices  in  our test suite}
  \label{fp1} 
\end{figure}

Our test suite consists of the following  transition matrices: 
\begin{enumerate}
    \item\label{it1} The  transition rate matrix associated   with  the queuing model described in \cite{Dudin}. This is a complex queuing model, a
BMAP/PHF/1/N model with retrial system with finite buffer and non-persistent customers. We do not describe in detail the construction of this matrix, as it would take
some space, but refer the reader to [6, Sections 4.3 and 4.5]. The buffer size is denoted as $n$. The only change with
respect to the paper is that we fix the orbit size to a finite capacity $k$ (when the orbit
is full, customers leave the queue forever). We set $n=15$, which results in a block
upper Hessenberg matrix Q of size $N=110400$ with  $k\times k$ blocks of size $138$. 
\item\label{it2}  The transition rate matrix for the  model described in Example 1 of \cite{PSS}.  The model describes a time-sharing system with $n$ terminals which share the same computing resource.  The matrices 
are nearly completely decomposable  (NCD)  so that classical stationary iterative methods do not perform satisfactorily as the  spectral gap is 
pathologically close to unity. We set $n=80$ which gives a matrix of size $N=91881$.
\item\label{it3} The transition rate matrix for the  model described in Example 3 of \cite{PSS}.  The model describes a multi-class, finite buffer, priority system.  The buffer size is denoted as $n$. This model
can be applied to telecommunications modeling, and has been used to model ATM queueing
networks as discussed in \cite{SWJ,PSS}.  We note that the model parameters can be  selected so that the resulting
Markov chain is nearly completely decomposable.  We set $n=100$ which results in a matrix of size $N=79220$.
\item\label{it4} The transition rate matrix generated by  the set of mutual-exclusion  problems 
considered in \cite{FPS}. In these problems, $n$  distinguishable processes share a certain resource. Each of these processes alternates between a
sleeping state and a resource using state. However, the number of processes that
may concurrently use the resource is limited to $r$ where $1\leq r\leq n$ so that when
a process wishing to move from the sleeping state to the resource using state
finds $r$ processes already using the resource, that process fails to access the
resource and returns to the sleeping state.  We set $n=16$  and $r=12$ so that the transition matrix has size $N=64839$.
\end{enumerate}
All the  considered transition matrices  are sparse matrices.  In  Figure \ref{fp1} we show the spy plots of the matrices generated in  tests 1-4.  
The matrices are stored using the  compressed sparse column format.  With this method, only nonzero entries are kept in memory. However, despite evident merits 
this solution has also  some drawbacks.  In particular, we notice that  certain operations such as {\tt submat} calls become relatively expensive. 
Clearly, the size of the block partitioning of  the matrix  seriously affects the performance of iterative methods. There is an extensive literature on this topic
(see for instance \cite{NS,KLEVANS199523} and the references given therein). 
We have not tried the
partition algorithms described in \cite{NS,KLEVANS199523} and 
in this paper we explore the use of blocks of equal size $n_i=\ell$, $1\leq i\leq n-1$. 
A  test for the change of two consecutive iterates  to be less than  a 
prescribed tolerance or the number of iterations to be greater than a given bound  is used as  stopping criterion. In other words, we stop the iteration  if
\[
\parallel \B x^{(k)}-\B  x^{(k+1)}\parallel_1\leq \epsilon \ \vee \ k\geq{\rm maxit}.
\]
In all  the experiments reported below we have used  $\epsilon=-1.0e\!-9$ and maxit=$1.0e\!+4$. For each algorithm and experiment  we measure the sequential completion time 
$T_{seq}$, the parallel completion time on $m$ threads $T_{par}(m)$, the speedup
$S_p(m)=T_{seq}/T_{par}(m)$ and the efficiency $E(m)=S_p(m)/m$.  We also report a plot of the residual $\parallel \B x^{(k)}-\B  x^{(k+1)}\parallel_1$ to analyze the convergence of the iterative scheme. 

\begin{figure}[ht] 
  \begin{subfigure}[b]{0.5\linewidth}
    \centering
    \fontsize{6}{10}\fontfamily{ptm}\selectfont
    \includegraphics[width=0.95\linewidth]{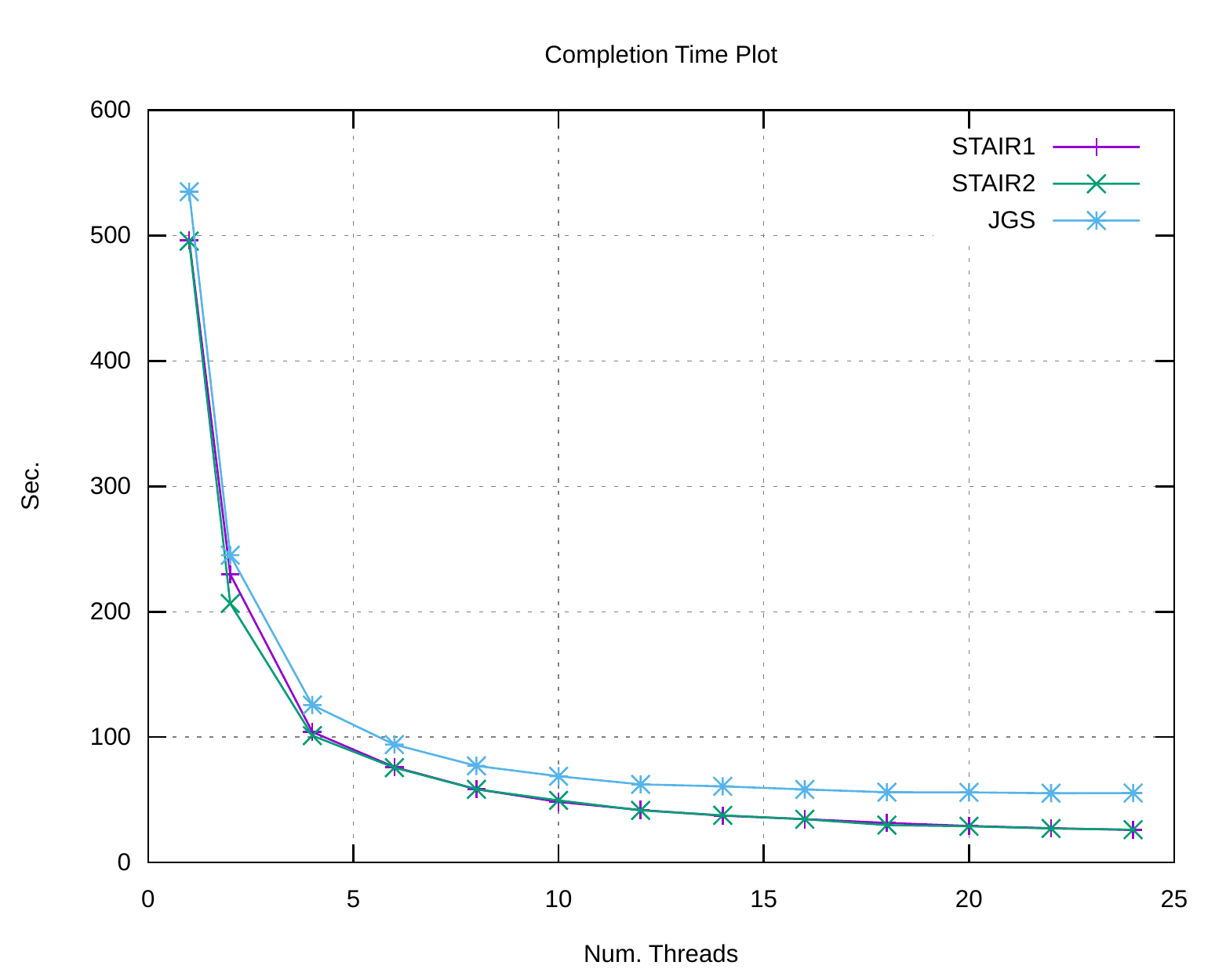}
    \caption{Completion-time plot for the matrix in \ref{it1}} 
    \label{fig9:a} 
    \vspace{4ex}
  \end{subfigure}
  \begin{subfigure}[b]{0.5\linewidth}
    \centering
    \fontsize{6}{10}\fontfamily{ptm}\selectfont
    \includegraphics[width=0.95\linewidth]{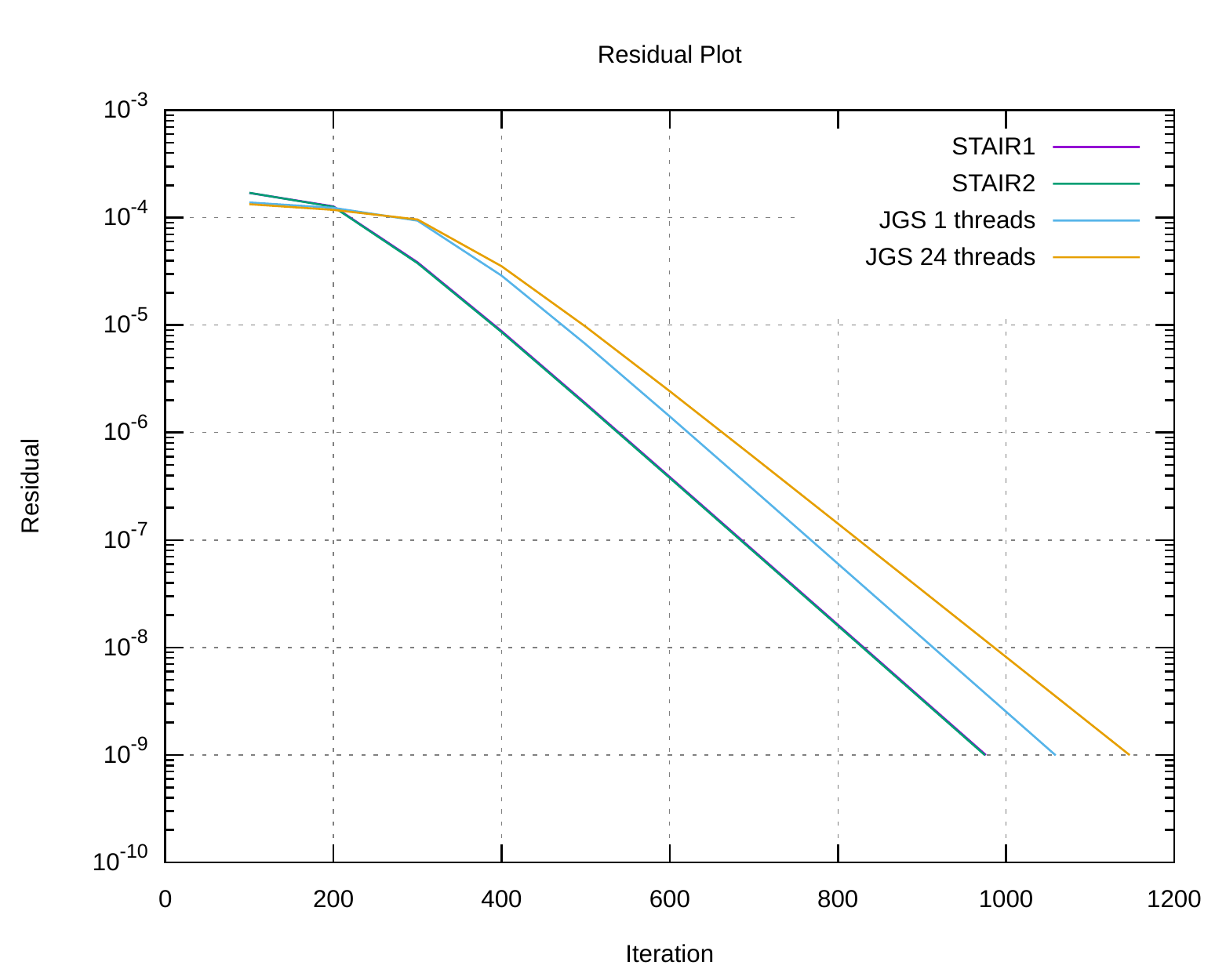}
    \caption{ Residual plot  for  the matrix in  \ref{it1}}
    \label{fig9:b} 
    \vspace{4ex}
  \end{subfigure} 
  \begin{subfigure}[b]{0.5\linewidth}
    \centering
    \fontsize{6}{10}\fontfamily{ptm}\selectfont
    \includegraphics[width=0.95\linewidth]{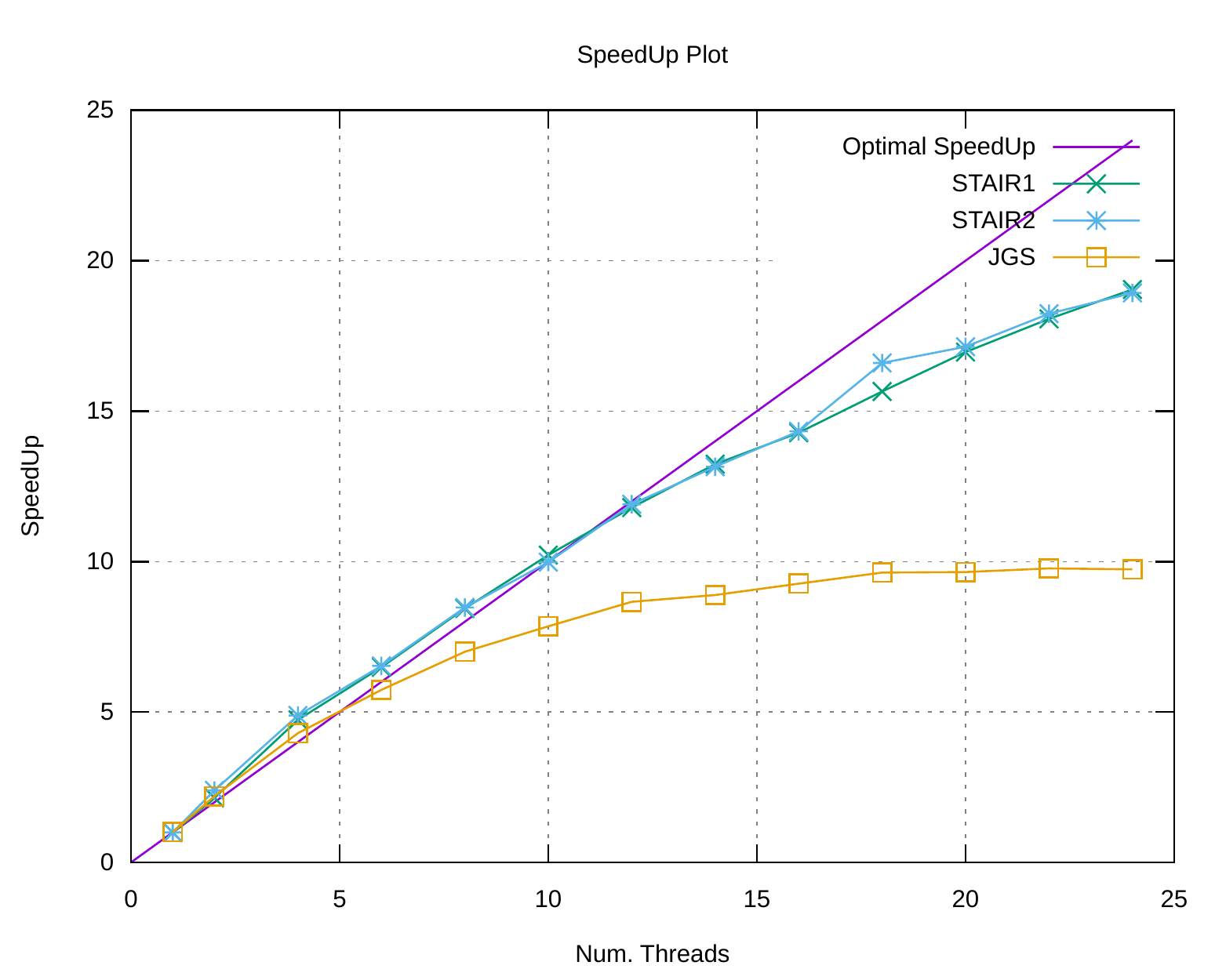}
      \caption{Speedup plot for  the matrix in  \ref{it1}} 
    \label{fig9:c} 
  \end{subfigure}
  \begin{subfigure}[b]{0.5\linewidth}
    \centering
    \fontsize{6}{10}\fontfamily{ptm}\selectfont
    \includegraphics[width=0.95\linewidth]{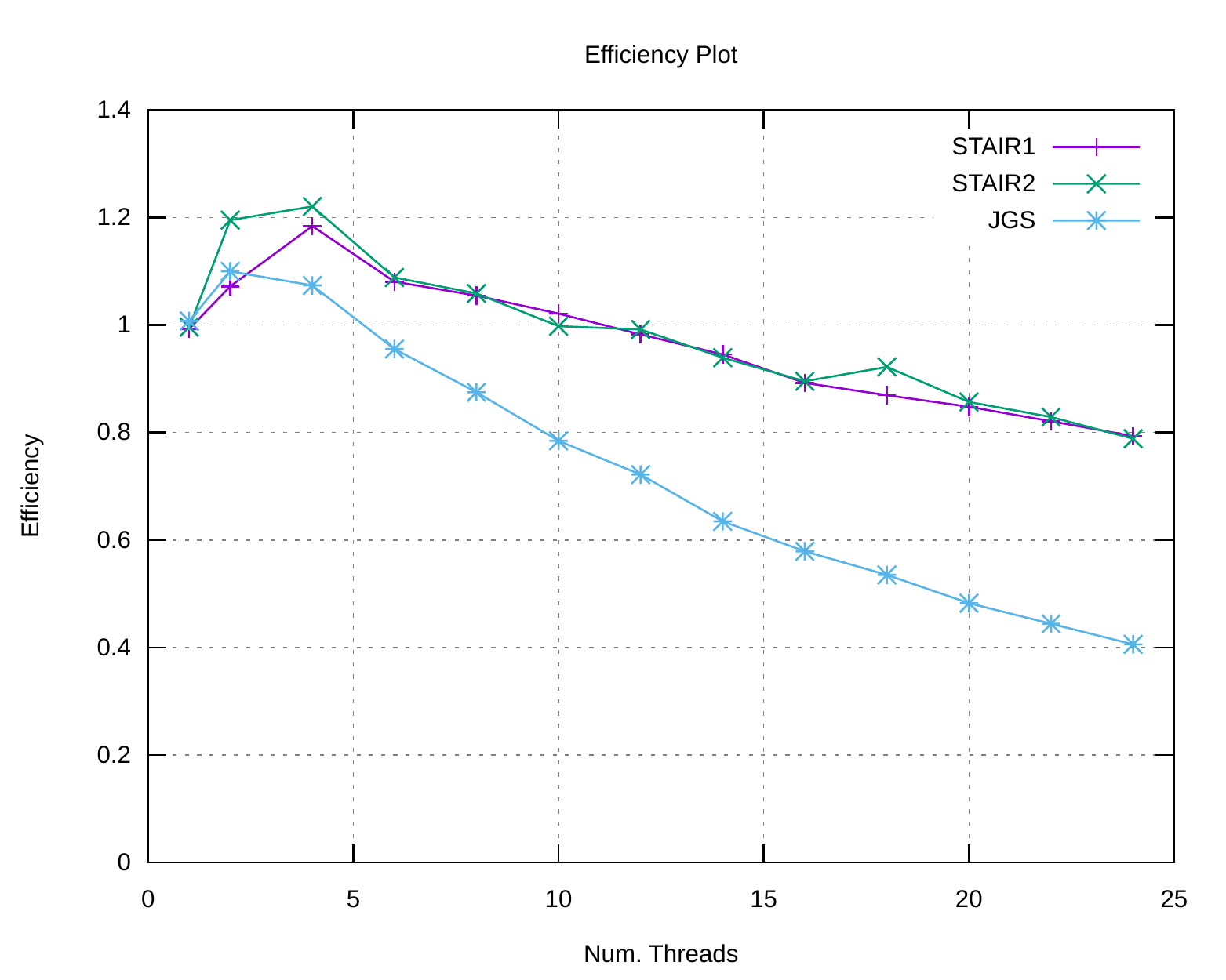}
    \caption{Efficiency  plot for  the matrix in  \ref{it1}}
    \label{fig9:d} 
  \end{subfigure} 
  \caption{Illustration of the  performance of algorithms {\tt JGS}, {\tt STAIR1} and {\tt STAIR2}  for the matrix in \ref{it1}}
  \label{fp3} 
\end{figure}

In the first experiment we consider the  transition rate matrix generated in \ref{it1} with $n=15$ and $k=800$.  The matrix has size $N=110400$  and the block partitioning is determined by  setting $\ell=512$.  The matrix is block lower Hessenberg and block banded.  In Figure \ref{fp3} we show  the plots  of completion time, residual, speedup and efficiency  generated  for this matrix.  The convergence of block staircase methods  is better than the convergence of the block Gauss-Seidel method.

\begin{figure}[ht] 
  \begin{subfigure}[b]{0.5\linewidth}
    \centering
    \fontsize{6}{10}\fontfamily{ptm}\selectfont
     \includegraphics[width=0.95\linewidth]{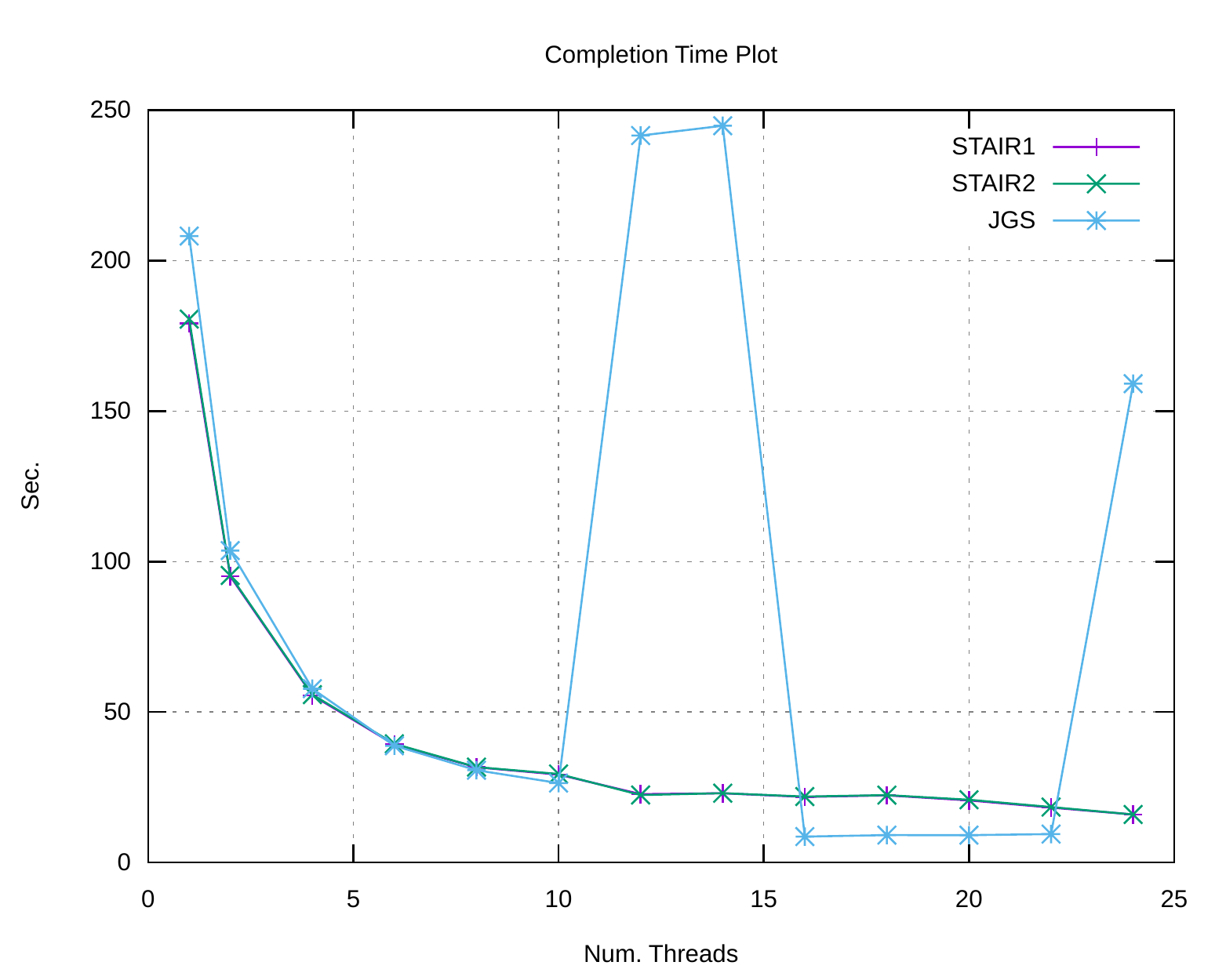}
    \caption{Completion-time  plot for the matrix in  \ref{it2}} 
    \label{fig10:a} 
    \vspace{4ex}
  \end{subfigure}
  \begin{subfigure}[b]{0.5\linewidth}
    \centering
    \fontsize{6}{10}\fontfamily{ptm}\selectfont
     \includegraphics[width=0.95\linewidth]{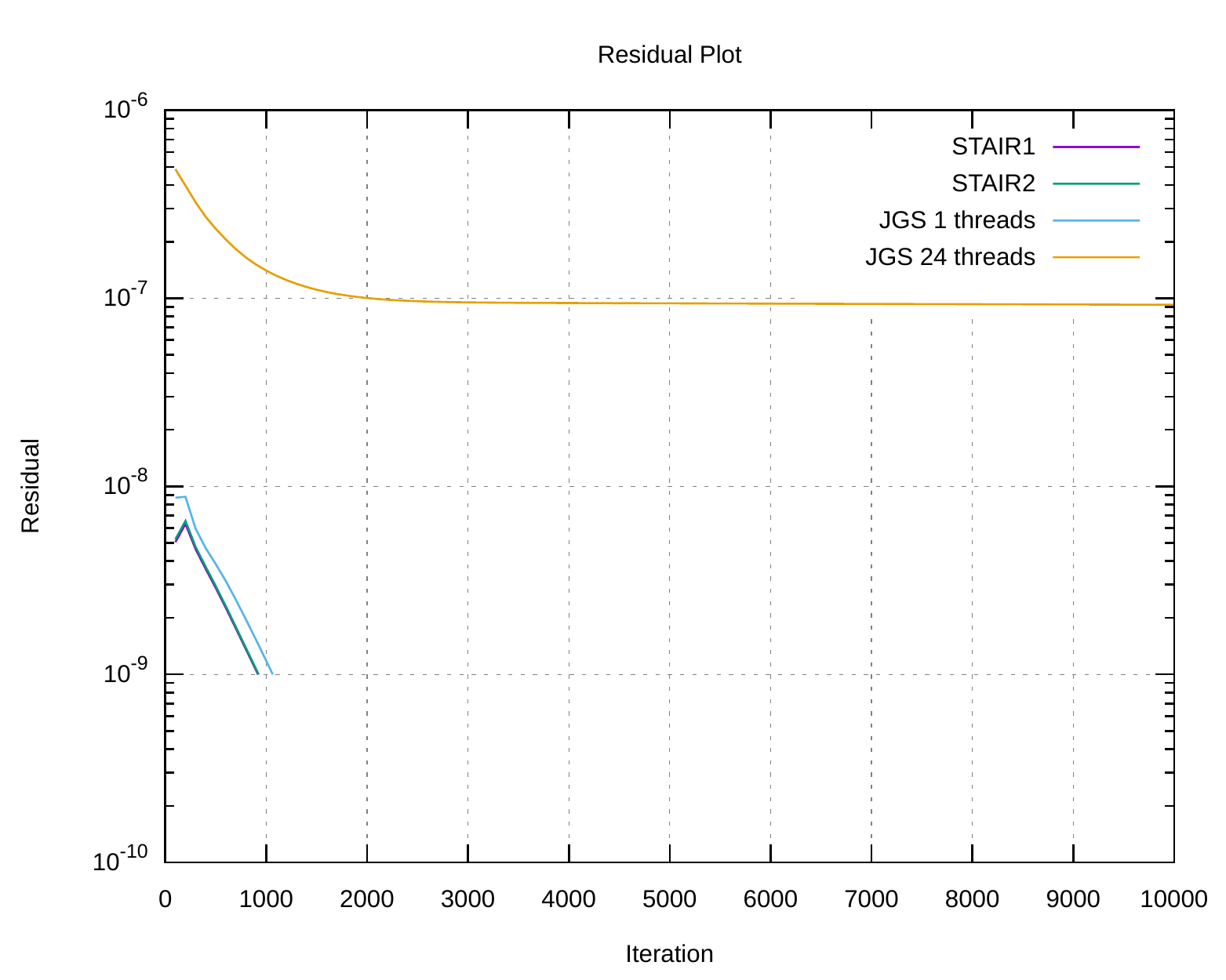}
    \caption{ Residual plot  for  the matrix in  \ref{it2}}
    \label{fig10:b} 
    \vspace{4ex}
  \end{subfigure} 
  \begin{subfigure}[b]{0.5\linewidth}
    \centering
    \fontsize{6}{10}\fontfamily{ptm}\selectfont
     \includegraphics[width=0.95\linewidth]{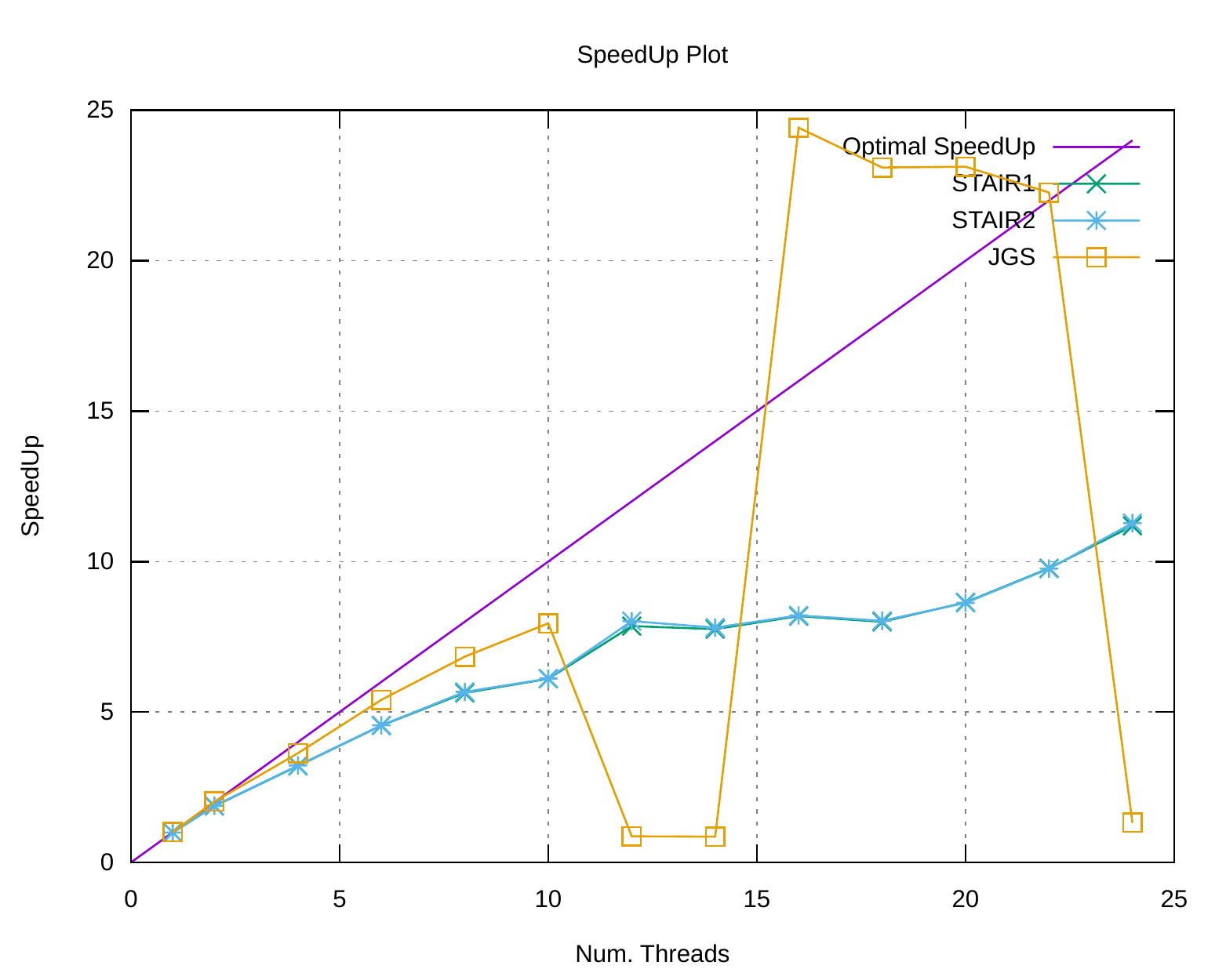}
      \caption{Speedup plot for  the matrix in  \ref{it2}} 
    \label{fig10:c} 
  \end{subfigure}
  \begin{subfigure}[b]{0.5\linewidth}
    \centering
    \fontsize{6}{10}\fontfamily{ptm}\selectfont
     \includegraphics[width=0.95\linewidth]{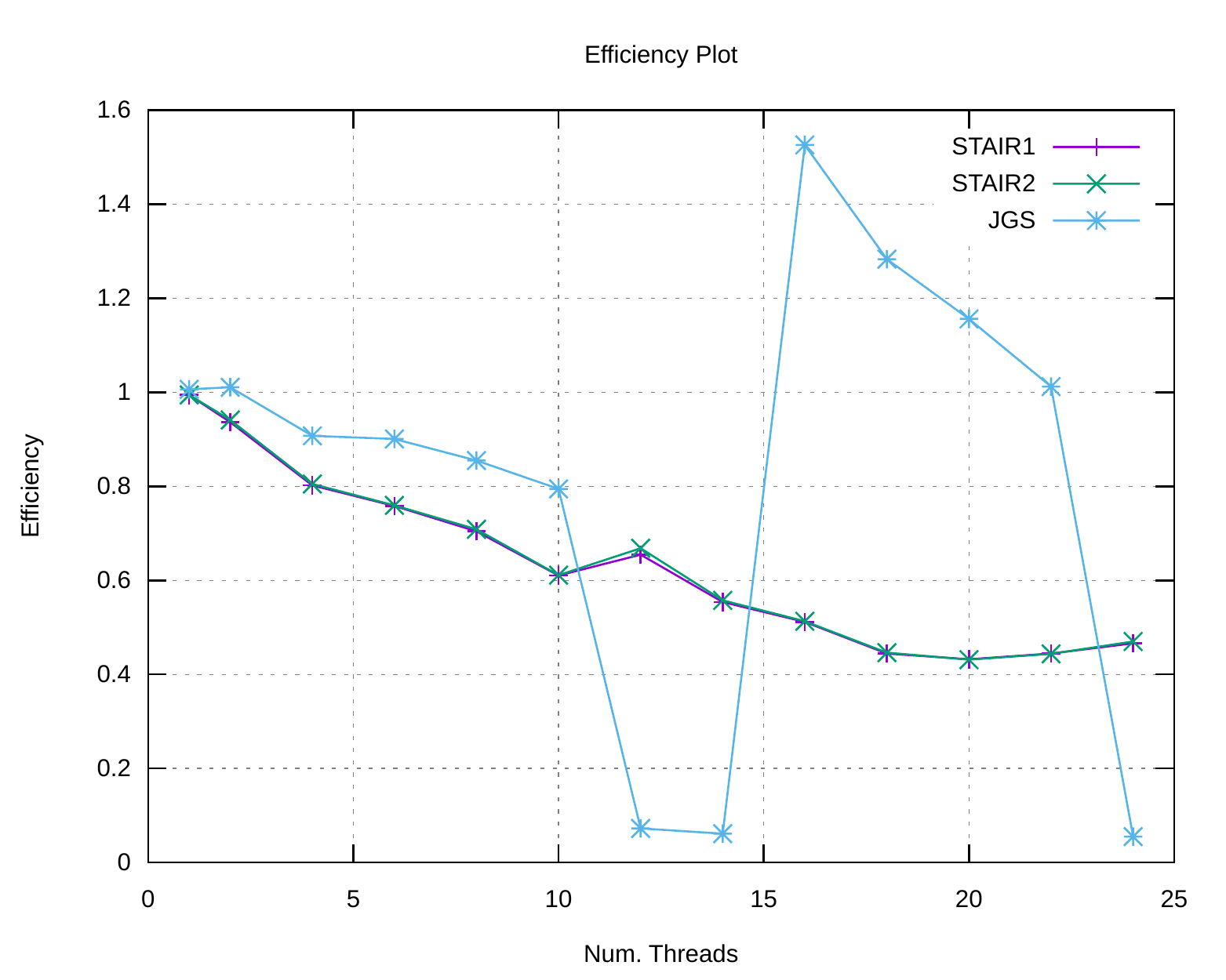}
    \caption{ Efficiency  plot for  the matrix in  \ref{it2}}
    \label{fig10:d} 
  \end{subfigure} 
  \caption{Illustration of   the  performance of algorithms {\tt JGS}, {\tt STAIR1} and {\tt STAIR2}  for the matrix in \ref{it2}}
  \label{fp4} 
\end{figure}

In Figure \ref{fp4} we show  the plots  generated  for the matrix in \ref{it2} with $n=80$.  The size is 
$N=91881$. The matrix is symmetric with bandwidth 1135. We set $\ell=2048$ so that the matrix is  block tridiagonal. 

\begin{figure}[ht] 
  \begin{subfigure}[b]{0.5\linewidth}
    \centering
    \fontsize{6}{10}\fontfamily{ptm}\selectfont
    \includegraphics[width=0.95\linewidth]{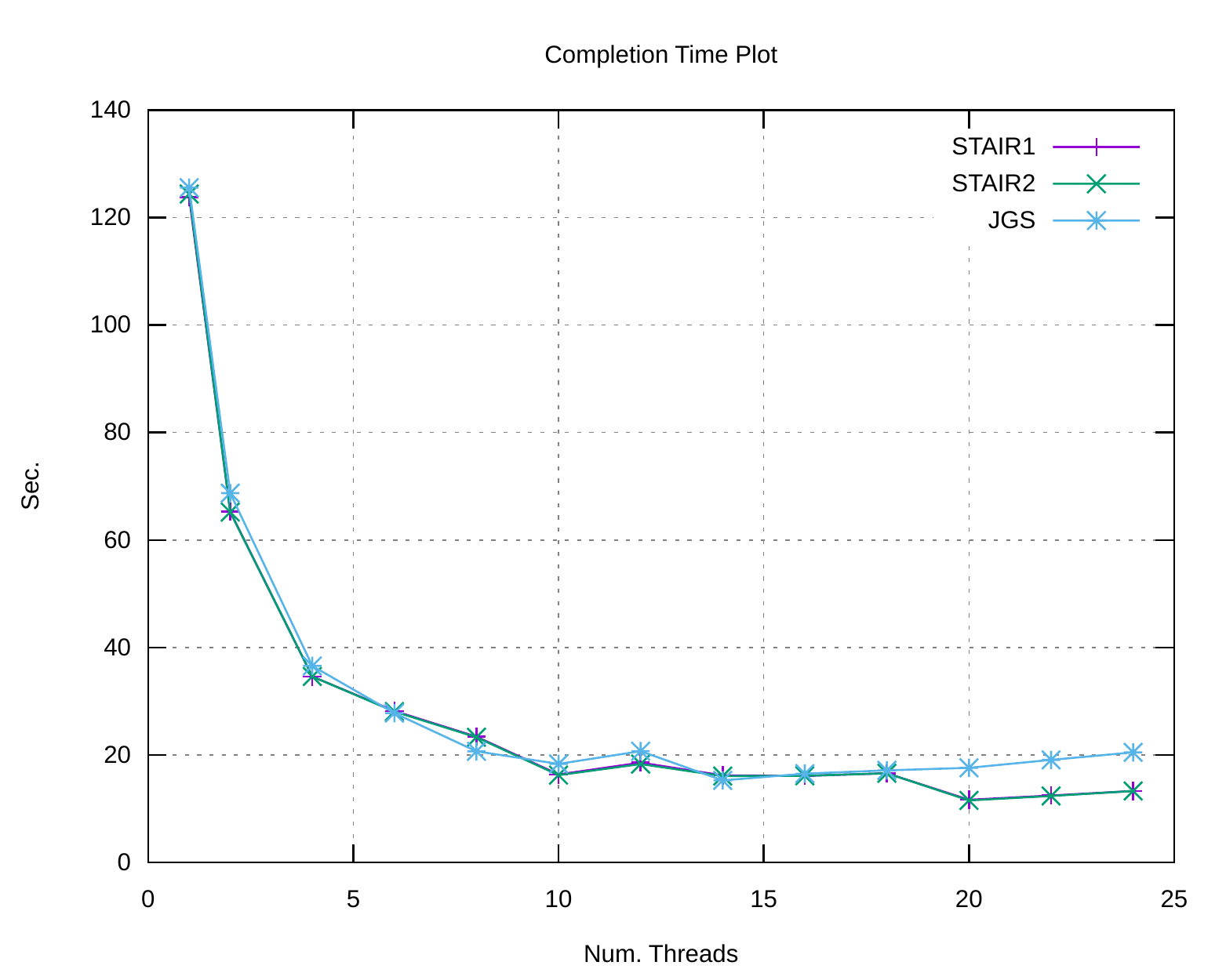}
    \caption{Completion-time plot for the matrix in \ref{it3}} 
    \label{fig11:a} 
    \vspace{4ex}
  \end{subfigure}
  \begin{subfigure}[b]{0.5\linewidth}
    \centering
    \fontsize{6}{10}\fontfamily{ptm}\selectfont
 \includegraphics[width=0.95\linewidth]{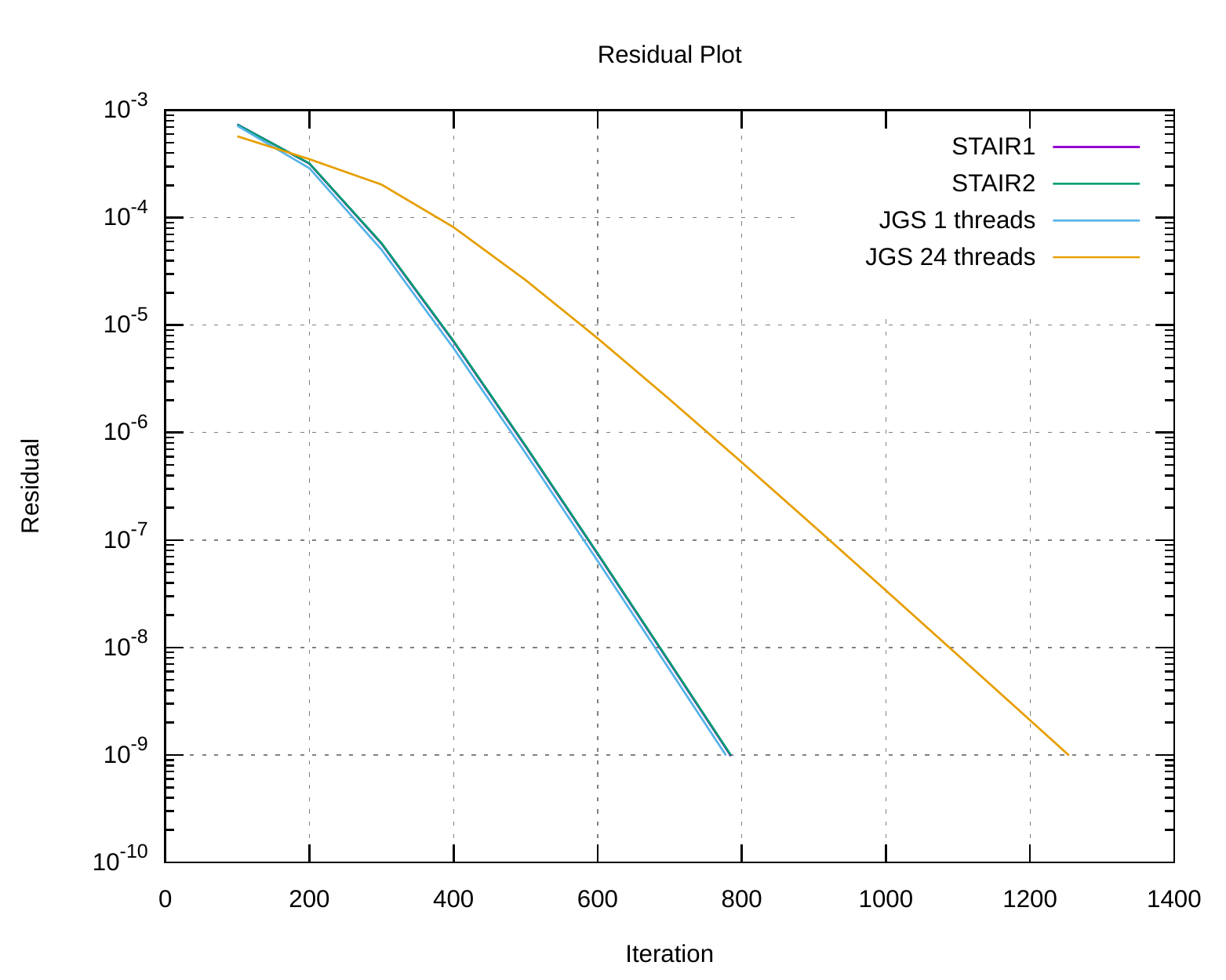}
    \caption{ Residual plot  for  the matrix in  \ref{it3}}
    \label{fig11:b} 
    \vspace{4ex}
  \end{subfigure} 
  \begin{subfigure}[b]{0.5\linewidth}
    \centering
    \fontsize{6}{10}\fontfamily{ptm}\selectfont
 \includegraphics[width=0.95\linewidth]{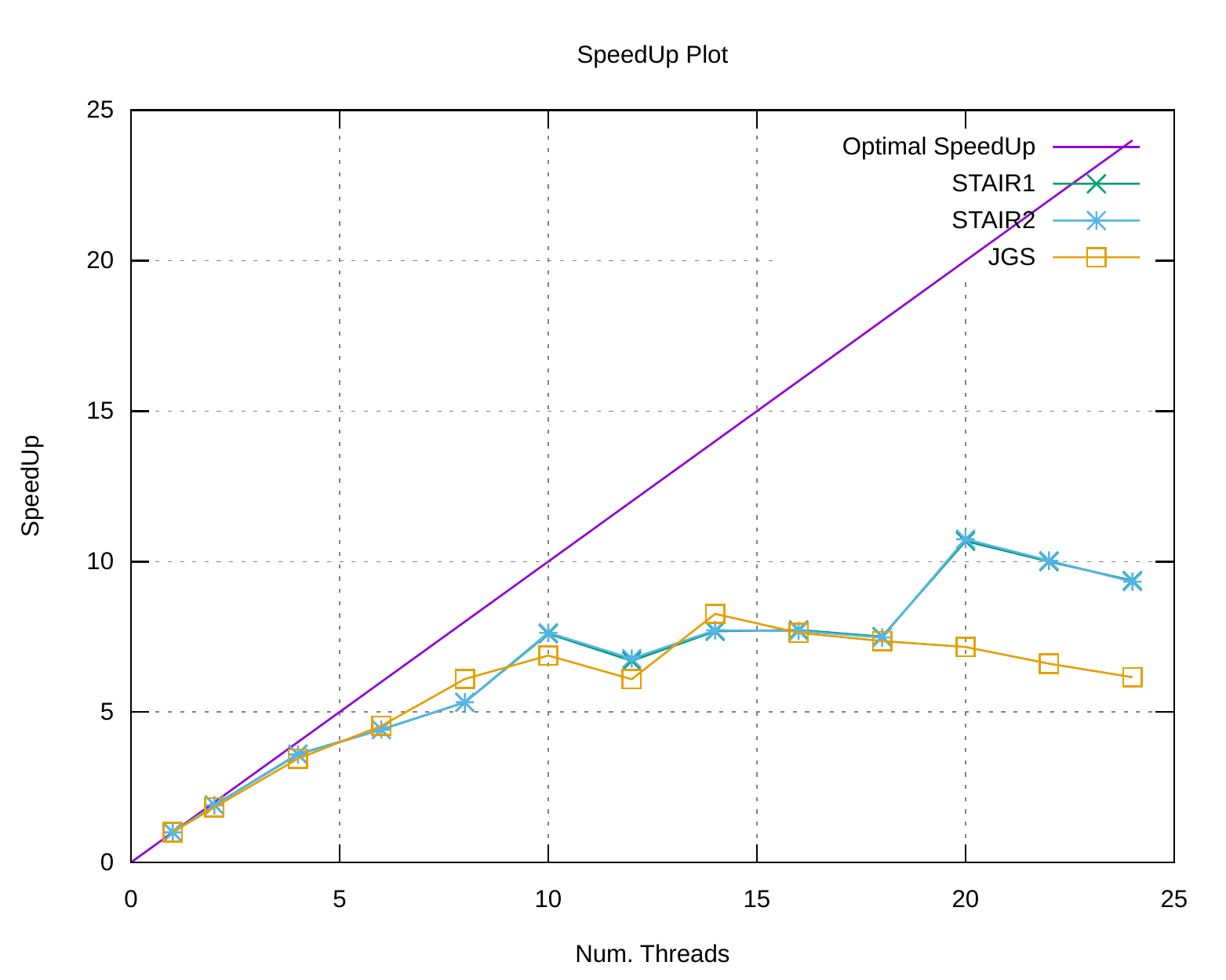}
    \caption{Speedup plot for  the matrix in  \ref{it3}} 
    \label{fig11:c} 
  \end{subfigure}
  \begin{subfigure}[b]{0.5\linewidth}
    \centering
    \fontsize{6}{10}\fontfamily{ptm}\selectfont
     \includegraphics[width=0.95\linewidth]{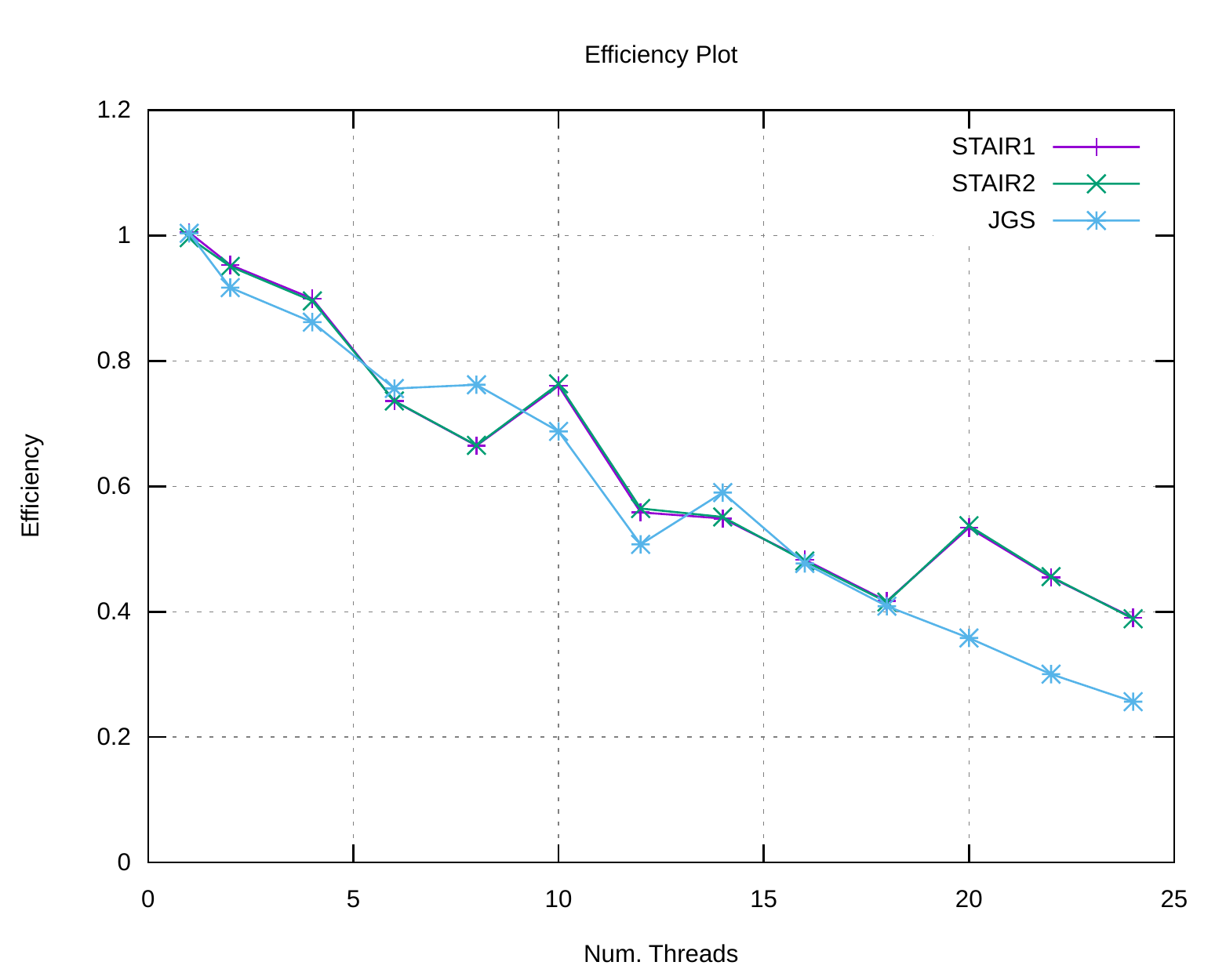}
    \caption{ Efficiency  plot for  the matrix in  \ref{it3}}
    \label{fig11:d} 
  \end{subfigure} 
  \caption{Illustration of   the  performance of algorithms {\tt JGS}, {\tt STAIR1} and {\tt STAIR2}  for the matrix in \ref{it3}}
  \label{fp5} 
\end{figure}

In Figure \ref{fp5} we show  the plots  generated  for the matrix in \ref{it3} with $n=100$. The matrix has size $N=79220$, lower bandwidth 2370 and upper bandwidth 1585. We set $\ell=2048$ so that the matrix is block banded in block  lower Hessenberg form. 

The results in Figure \ref{fp4}, \ref{fp5}  clearly highlight that the convergence of the 
{\tt JGS} method  may deteriorate as the number of threads increase since the iteration becomes close to a pure block Jacobi method.  In the case \ref{it2} of a nearly completely decomposable  transition rate matrix the deterioration can produce divergence phenomena. In particular for the number of threads $m\in {12, 14, 24}$  {\tt JGS} does not converge and the process  stop as the   maximum number of iterations has been reached.

\begin{figure}[ht] 
  \begin{subfigure}[b]{0.5\linewidth}
    \centering
    \fontsize{6}{10}\fontfamily{ptm}\selectfont
     \includegraphics[width=0.95\linewidth]{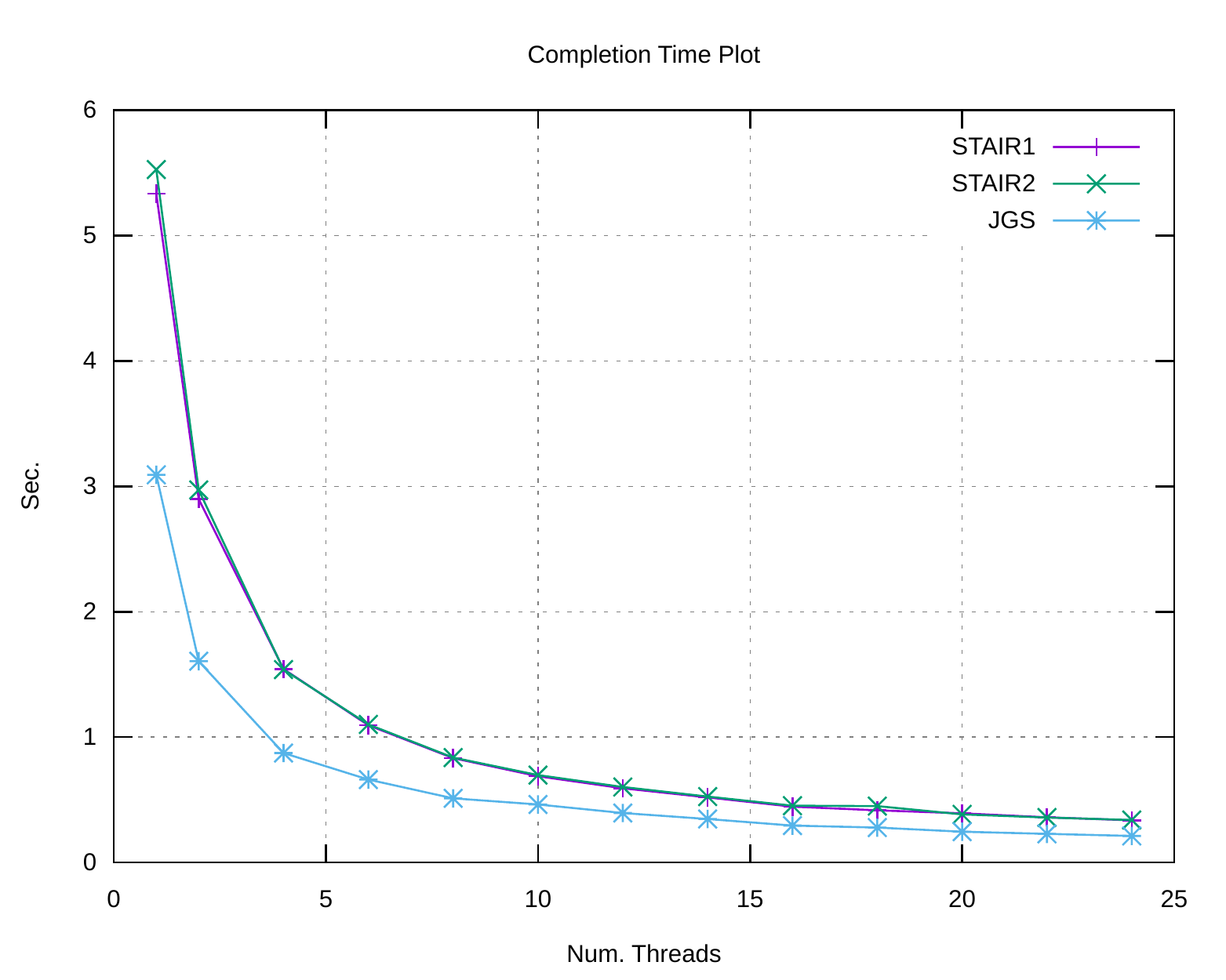}
    \caption{Completion-time  plot for the matrix in  \ref{it4}} 
    \label{fig12:a} 
    \vspace{4ex}
  \end{subfigure}
  \begin{subfigure}[b]{0.5\linewidth}
    \centering
    \fontsize{6}{10}\fontfamily{ptm}\selectfont
     \includegraphics[width=0.95\linewidth]{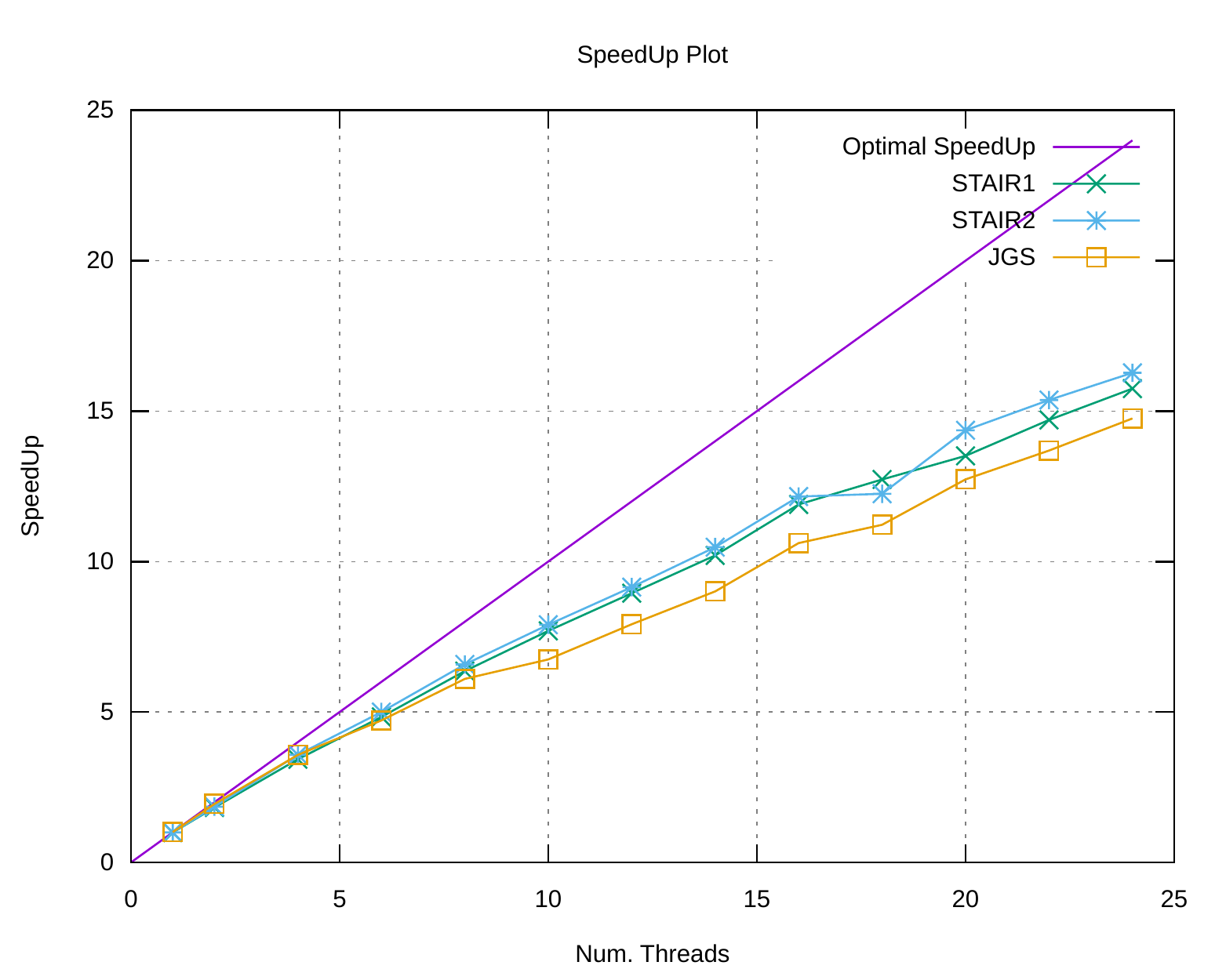}
    \caption{ Speedup plot  for  the matrix in  \ref{it4}}
    \label{fig12:b} 
    \vspace{4ex}
  \end{subfigure} 
  \caption{Illustration of   the  performance of algorithms {\tt JGS}, {\tt STAIR1} and {\tt STAIR2}  for the matrix in \ref{it4}}
  \label{fp6} 
\end{figure}

Finally, in Figure \ref{fp6} we illustrate  the plots  of completion time and speedup for the matrix generated  by \ref{it4} with $n=16$, and $r=12$. The matrix has size $N=64839$ and bandwidth 13495. We set $\ell=256$ so that the matrix is block banded.  Since the entries are very rapidly decaying away from the main diagonal  all methods perform quite well and the number of iterations in the {\tt JGS} method is quite insensitive to the number of threads. 

\begin{figure}[ht]
  \begin{subfigure}[b]{0.5\linewidth}
    \centering
    \fontsize{6}{10}\fontfamily{ptm}\selectfont
    \includegraphics[width=0.95\linewidth]{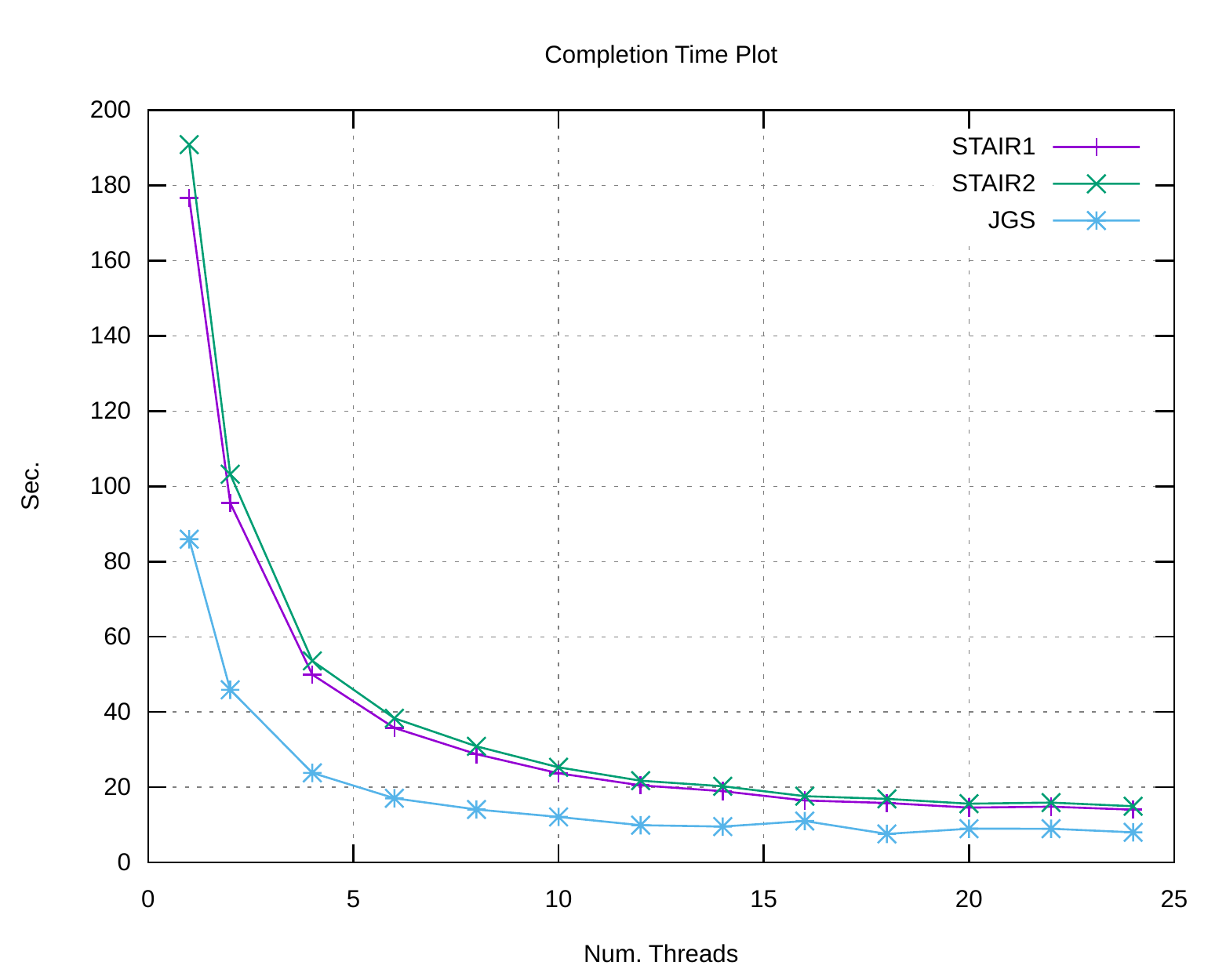}
    \caption{Completion-time plot  with $\ell=512$}
    \label{fig13:a} 
    \vspace{4ex}
  \end{subfigure}
  \begin{subfigure}[b]{0.5\linewidth}
    \centering
    \fontsize{6}{10}\fontfamily{ptm}\selectfont
     \includegraphics[width=0.95\linewidth]{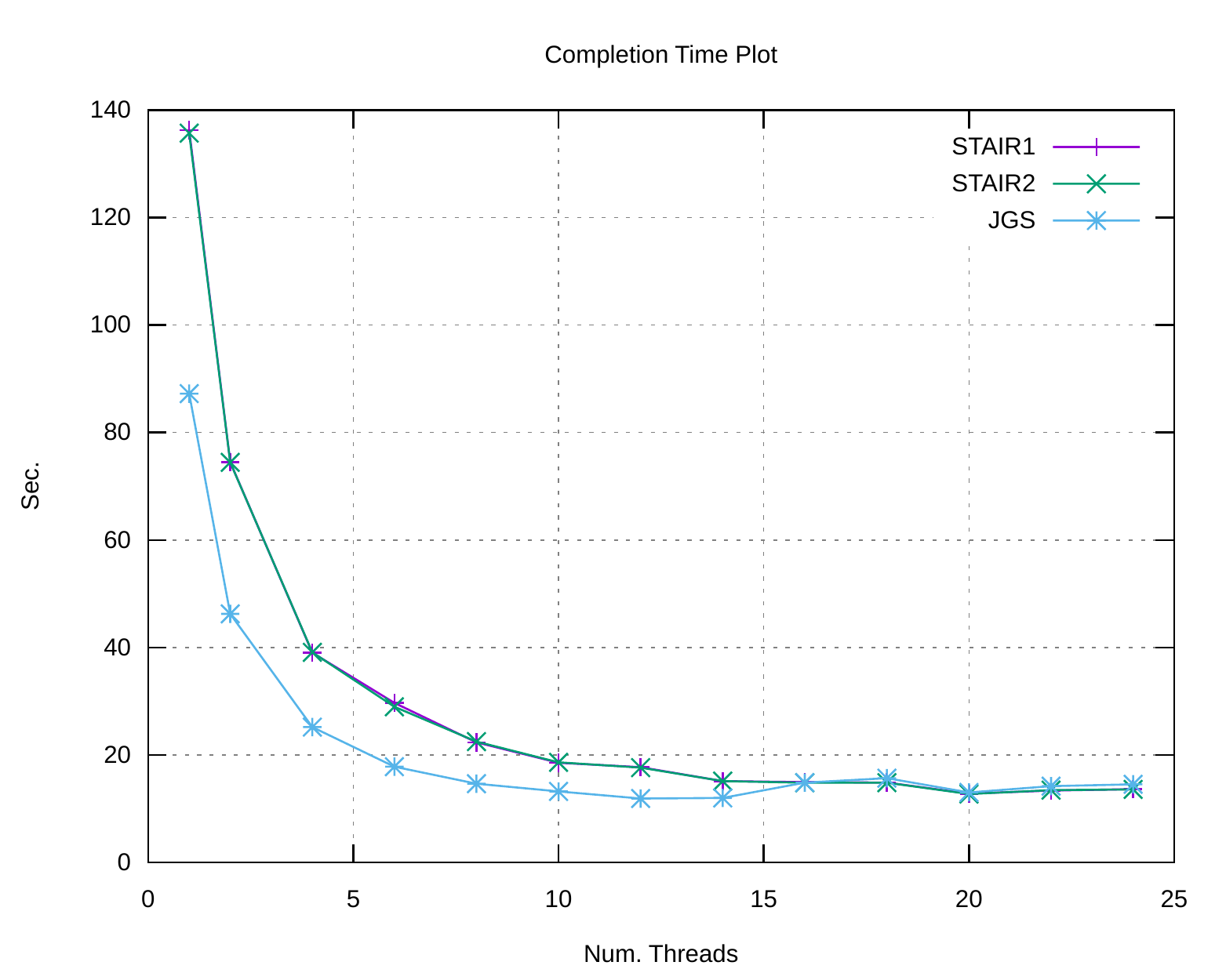}  
    \caption{Completion-time plot  with $\ell=1024$}
    \label{fig13:b} 
    \vspace{4ex}
  \end{subfigure}
  \begin{subfigure}[b]{0.5\linewidth}
    \centering
    \fontsize{6}{10}\fontfamily{ptm}\selectfont
     \includegraphics[width=0.95\linewidth]{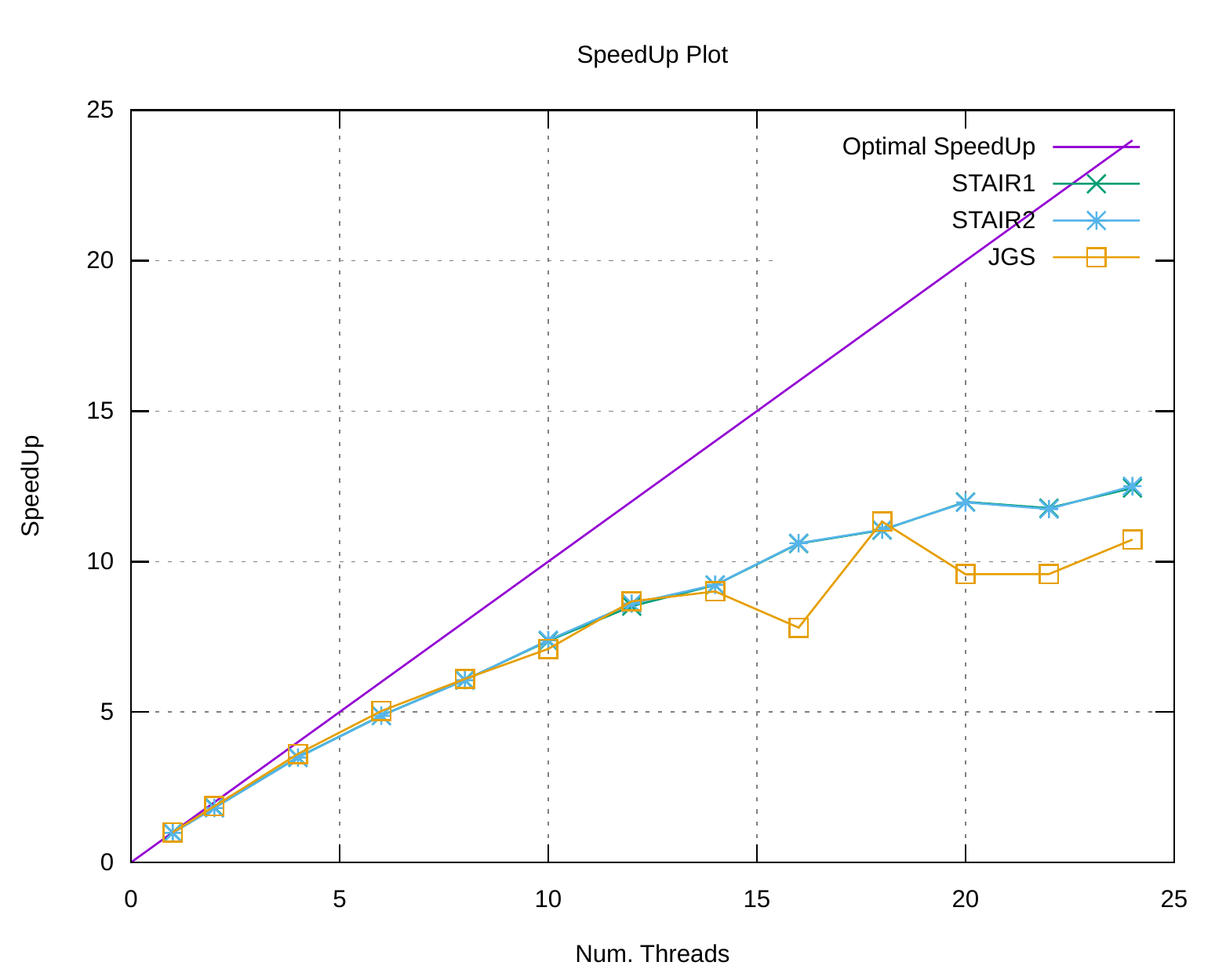}
    \caption{Speedup plot with $\ell=512$}
    \label{fig14:a} 
    \vspace{4ex}
  \end{subfigure}
  \begin{subfigure}[b]{0.5\linewidth}
    \centering
    \fontsize{6}{10}\fontfamily{ptm}\selectfont
     \includegraphics[width=0.95\linewidth]{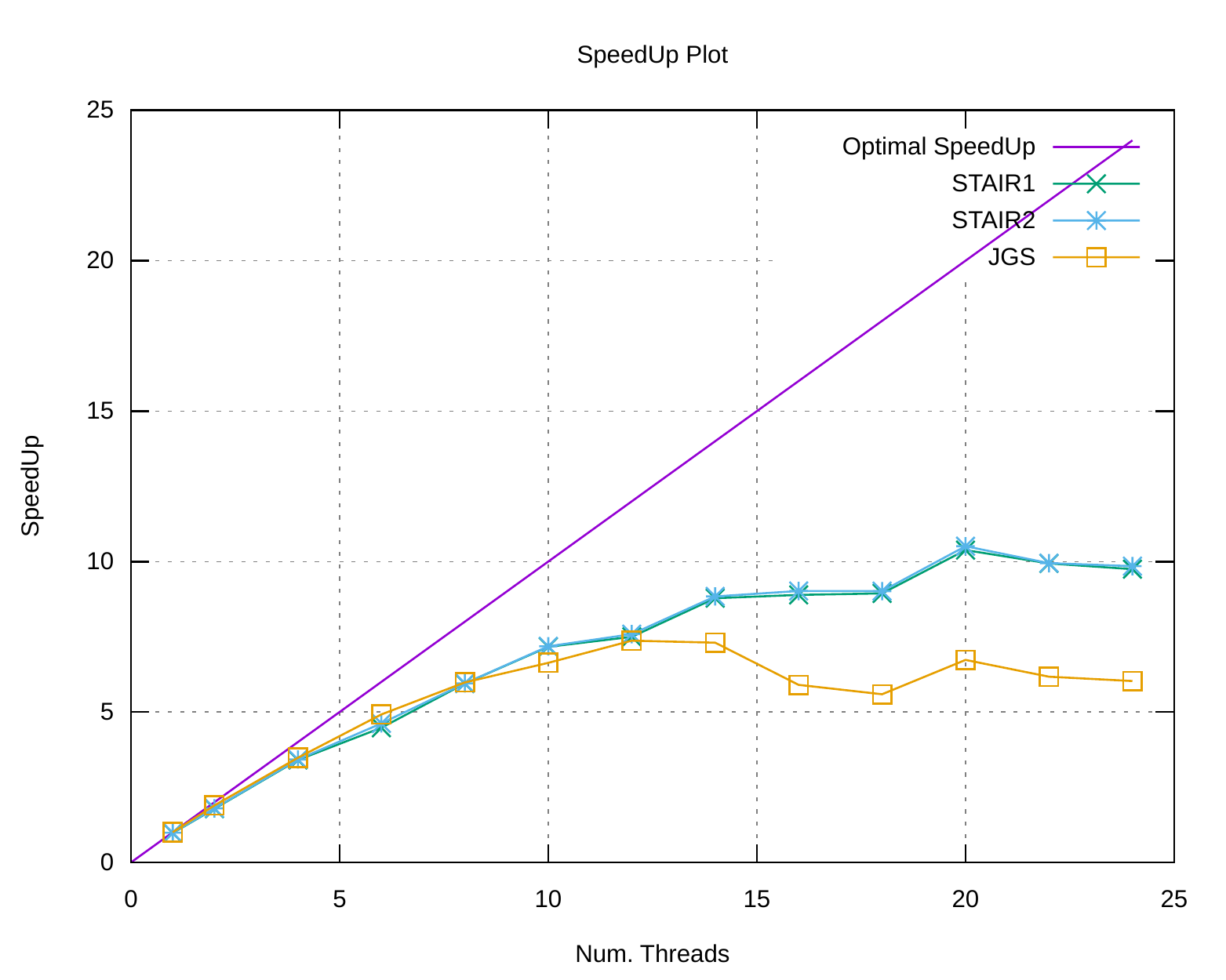}
    \caption{Speedup plot  with $\ell=1024$}
    \label{fig14:b} 
    \vspace{4ex}
  \end{subfigure} 
  \caption{Illustration of the completion time and speedup plot for the matrix \ref{it3} with different block sizes}
  \label{fp7} 
\end{figure}

Concerning the parallel  performance, we 
recall  that the server has only  24 combined physical cores, and going above 12  required  communication  between  the  different  CPUs,  which  inevitably  reduces the efficiency of the parallelization.  When the number of threads is quite small  it is  generally observed that the bigger the block size, the shorter is  the execution time. Differently, as the number of threads increases small blocks  promote  the parallelism.    In Figure \ref{fp7} we show the the speedup plot for the test \ref{it3} with $\ell=512$ and $\ell=1024$, respectively.  The comparison with the results  reported in Figure \ref{fp5}
 with $\ell =2048$ indicates some improvements.  These  effects are  enlightened  in all the conducted experiments.  Considering that most consumer hardware has between 2 and 8 or 16 cores, this shows that the  proposed method is  generally well tuned for the currently available architectures.

\section{Conclusions}
This paper presents the results of a preliminary experimental investigation of the  performance of a stationary iterative method based on a block staircase splitting  for solving singular systems of linear equations arising in Markov chain modelling.  From the experiments presented, we can deduce  that the 
method is well suited for solving block banded or more generally localized systems  in a shared-memory  parallel computing environment.  The parallel implementation has been benchmarked using several Markovian models. In the future we plan to examine the performance of block staircase splittings in a distributed computing environment and, moreover,  their use as preconditioners for other iterative methods.


\end{document}